\documentclass[times, doublespace]{cnmauth}

\usepackage[english]{babel}
\usepackage{moreverb}

\newcommand\BibTeX{{\rmfamily B\kern-.05em \textsc{i\kern-.025em b}\kern-.08em
T\kern-.1667em\lower.7ex\hbox{E}\kern-.125emX}}

\usepackage{hyperref}
\usepackage{float}
\usepackage{amssymb}
\usepackage{amsmath}
\usepackage{latexsym}
\usepackage{amsthm}
\usepackage{mathrsfs} 
\usepackage{fancyhdr}
\usepackage{graphicx}
\usepackage{newlfont}
\usepackage{textcomp}
\usepackage{caption}
\usepackage{booktabs}
\usepackage{multirow}

\usepackage{hyperref}
\usepackage{color}

\usepackage{etoolbox}
\AtBeginEnvironment{bmatrix}{\setlength{\arraycolsep}{1pt}} % to make block matrices more compact
\makeatletter
\renewcommand*\env@matrix[1][*\c@MaxMatrixCols c]{%
  \hskip -\arraycolsep
  \let\@ifnextchar\new@ifnextchar
  \array{#1}}
\makeatother

\newcommand*{\bm}[1]{\boldsymbol{#1}}
\newcommand{\aled}[1]{\widehat{#1}}

\begin{document}

\runningheads{S. Calandrini et al.}{Fluid-structure interaction simulations of venous valves}

\title{Fluid-structure interaction simulations of venous valves: a monolithic ALE method for large structural displacements}

\author{S. Calandrini \affil{1} \corrauth, E. Aulisa \affil{1} }

\address{\affilnum{1} Department of Mathematics and Statistics, Texas Tech University, Lubbock TX 79409, USA.}

\corraddr{Sara Calandrini, Broadway \& Boston, Department of Mathematics and Statistics, Texas Tech University, Lubbock TX 79409, USA.  E-mail: sara.calandrini@ttu.edu}

\cgsn{National Science Foundation}{DMS-1412796}

\begin{abstract} 
Venous valves are bicuspidal valves that ensure that blood in veins only flows back to the heart. 
To prevent retrograde blood flow, the two intraluminal leaflets meet in the center of the vein and occlude the vessel.
In fluid-structure interaction (FSI) simulations of venous valves, the large structural displacements may lead to mesh 
deteriorations and entanglements, causing instabilities of the solver and, consequently,
the numerical solution to diverge. In this paper, we propose an Arbitrary Lagrangian-Eulerian (ALE) scheme for FSI
simulations designed to solve these instabilities. 
A monolithic formulation for the FSI problem is considered and, due to
the complexity of the operators, the exact Jacobian matrix is evaluated using automatic differentiation.
The method relies on the introduction of a staggered in time velocity %in the discretization scheme 
to improve stability, and on fictitious springs to model the contact force of the valve leaflets.
Since the large structural displacements may compromise the quality of the fluid mesh as well, 
a smoother fluid displacement, obtained with the introduction of a scaling factor that measures the distance of a fluid element
from the valve leaflet tip, guarantees that there are no mesh entanglements in the fluid domain.
To further improve stability, a Streamline Upwind Petrov Galerkin (SUPG) method is employed.
The proposed ALE scheme is applied to a 2D model of a venous valve. The presented simulations show 
that the proposed method deals well with the large structural displacements of the problem, 
allowing a reconstruction of the valve behavior in both the opening and closing phase.  
% especially in the closing phase.
\end{abstract}
 

\maketitle

\section{Introduction} 
Blood vessels allow blood to flow quickly and efficiently inside the human body. 
The three major types of blood vessels are: arteries, capillaries and veins. 
While arteries carry highly oxygenated blood from the heart to the rest of the body, veins carry de-oxygenated blood back to the heart. 
% Capillaries are the smallest and thinnest of the blood vessels in the body, as well as the most common. Capillaries carry blood very % close to the cells of the tissues of the body in order to exchange gases, nutrients, and waste products.
% On the other hand, veins act as the blood return counterparts of arteries, namely they carry de-oxygenated blood back to the heart.
Arteries face high levels of blood pressure as they carry blood being pushed from the heart.
On the other hand, veins are subjected to very low blood pressures and rely on muscle contractions and one-way valves to push blood back to the heart.
The muscle pump effect, which denotes the mechanism of venous return via muscle contraction, has been intensively studied from a 
medical point of view (\cite{laughlin1987skeletal, rowell1993human, naadland2009effect}). Numerical studies of
this mechanism have also been performed in \cite{marchandise2010accurate, fullana2009branched, keijsers20151d}.
Valves play an important role during the muscle pump effect. Their main contribution is to prevent backflow.  
Venous valves are bicuspidal valves, meaning they are composed of two intraluminal membranes (or leaflets), that meet in the center of the vein and occlude the vessel to prevent backflow. 
The physics of vein valves has been studied by Lurie et al. in \cite{lurie2002mechanism, lurie2003mechanism} and \cite{lurie2006}. 
In the latter, it is postulated that a valve cycle is composed of four phases: the opening phase, the equilibrium phase, the closing phase and the closed phase. The equilibrium phase is a complex stage of a valve cycle, during which the valve is maximally open. 
The leaflets are at some distance from the wall and this creates a narrowing of the lumen where the flow accelerates resulting in a proximally directed flow jet. In the sinus pocket behind the valve cusps, blood stream forms a vortex along the sinus wall and
the mural side of the valve cusps. These vortexes create a rising pressure on the mural side and a falling pressure
on the luminal side of the cusps initiating the closing phase.

In fluid-structure interaction (FSI) simulations of venous valves, the complex interactions between blood, vessels walls and valve leaflets are not the only challenge. The large structural displacements may lead to mesh deterioration and entanglements, causing instabilities of the solver and the numerical solution to diverge.
To deal with large structural displacements, several approaches have been proposed in the literature, including
immersed boundary methods (\cite{peskin2002immersed, peskin1989three}), Lagrange multiplier based fictitious domain methods (\cite{baaijens2001fictitious, de2003three, van2004combined, van2005three}), 
Lattice-Boltzmann methods (\cite{krafczyk1998analysis, fang2002lattice, feng2004immersed, buxton2006computational}) and standard ALE algorithms with remeshing (\cite{margolin2003second, loubere2010reale}). 
Many of these works (\cite{krafczyk1998analysis, de2003three, van2004combined, van2005three}) involve the modeling of heart valves.
A review of the state-of-the-art for FSI methods applied to aortic valves simulations can be found in
\cite{marom2015numerical}.
Recently, an extended ALE method has been proposed in \cite{basting2017extended}, and it has been applied to a 2D FSI benchmark problem modeling valves. The main characteristic of such a method is a variational mesh optimization approach 
that does not rely on any combinatorial considerations.
Focusing on venous valves simulations, the first FSI study of vein valves reported in the literature is \cite{buxton2006computational}. In the aforementioned paper, the fluid is captured using a Lattice Boltzmann (LB) model, while the solid mechanics is captured using a Lattice Spring (LS) model. The dynamics of the valve opening area and the blood flow rate through the valve are investigated to capture the qualitative behavior of vein valves.
In \cite{narracott2010analysis}, both the dynamics of a mechanical heart valve prosthesis and the function of a native venous valve are analyzed. Two commercial codes are employed to study the distinct cardiovascular problems. The LS-DYNA explicit dynamics code is used for the venous valve simulations, where dependence of the venous valve function on the parent vessel elastic modulus is shown.
More recently, FSI simulations of venous valves have been performed in \cite{tien2014characterization, narracott2015fluid}.

In this paper, we focus on the large structural displacements that arise in FSI simulations of venous valves, proposing a monolithic ALE scheme that aims to solve the numerical instabilities caused by these large displacements.
The proposed method relies on three main features. First, we employ a staggered in time velocity in the time integration scheme, which improves the computational stability of the simulations.
This staggered approach is used in the momentum balance to compute the mesh velocity and in the kinematic equation to compute the velocity in the solid region.
Second, fictitious springs are used to model the contact force of the valve leaflets. Fictitious elastic springs 
are widely used in structural dynamics for the resolution of dynamics contact problems (\cite{gambhir2013stability, filiatrault2013elements}).
In a finite element setting, the contact of two structures can cause serious instabilities, and fictitious elastic springs allow the valve leaflets to get very close to each other without actually touching. 
Finally, since the large structural displacements may compromise the quality of the fluid mesh as well, a scaling factor that measures the distance between a fluid element and the valve leaflets is introduced. This scaling factor guarantees that there are no mesh entanglements in the fluid domain. To further improve stability, a Streamline Upwind Petrov Galerkin (SUPG) stabilization is employed. Following \cite{franca1993element}, the stabilization parameter, $\tau$, is computed solving a generalized eigenvalue problem.
The proposed method is applied to a 2D model of a venous valve. Two tests are performed, with the first one testing the mesh independence of the ALE scheme, and the second one assessing the sensitivity of valve dynamics to the elastic modulus of
the valve leaflet, independent of vein wall properties. 

The paper is structured as follows. In section \ref{formulation} the weak monolithic formulation of the FSI problem is described. In section \ref{supg} the SUPG stabilization added to the formulation presented in \ref{formulation} is illustrated. In section \ref{staggered} the time discretization scheme, which employes a staggered in time velocity, is described. In section \ref{functionK} the distance function used to handle possible entanglements in the fluid mesh is analyzed, while in section \ref{springs} fictitious springs are described, namely illustrating how springs are introduced in the formulation and the elastic parameters adopted.
In section \ref{tests}, the numerical tests are presented. Finally, in section \ref{fine} we draw our conclusions. 

\section{FSI Formulation}\label{formulation}
This section describes the FSI formulation of the venous valve problem. From now on, by solid structure we refer to the vessel walls and the valve leaflets, while the fluid considered is blood.  The solid motion is described in a
Lagrangian way, while the fluid is observed in Eulerian fashion. The deformation of the fluid domain
is taken into account according to an Arbitrary Lagrangian Eulerian (ALE) approach (\cite{FSIarticle, donea1982arbitrary, sackinger1996newton, wick2011fluid}). 
We consider a monolithic coupling between the fluid and the solid, focusing on a time-dependent formulation.
% An FSI problem can be formulated in terms of three subproblems: the fluid
% subproblem, the solid subproblem, and the subproblem for the fluid domain displacement.
To solve the FSI system, we use a monolithic Newton-Krylov solver preconditioned by a geometric multigrid algorithm (\cite{FSIarticle}). 
A validation of the FSI model and solver throughout a series of 2D and 3D benchmark tests can be found in \cite{FSIarticle},
where the results obtained were compared with several existing FSI models and solvers (\cite{turek2010numerical_2, richter2015monolithic}). Previous applications of our FSI model and
solver to hemodynamics problems can be found in \cite{COMPDYN2017, COUPLED2017} and \cite{calandrini2017magnetic}. 
Before discussing the formulation of the problem, we introduce some notations.
\subsection{Notation}
Let $\widehat{\Omega} = \widehat{\Omega}^f \cup \widehat{\Omega}^s \subset \mathbb{R}^n$ 
be a reference configuration consisting of a vein containing a bicuspidal valve,
where the superscripts $f$ and $s$ refer to fluid and solid, respectively.
Let $\Omega_t = \Omega^f_t \cup \Omega^s_t \subset \mathbb{R}^n$ be the current configuration at time $t$.
% Any reference configuration can be equivalently chosen in principle, and the natural
% choice is the initial configuration, $\Omega_0=\widehat{\Omega}$.
Let $\Gamma^i_t = \Omega^f_t \,\cap\, \Omega^s_t$ and $\widehat{\Gamma}^i = \widehat{\Omega}^f \,\cap\, \widehat{\Omega}^s$
be the interfaces between solid and fluid in the current and reference configuration, respectively.
Moreover, we define the parts of the boundary adjacent only to the fluid or only to the solid as
$\Gamma^f_t$, $\Gamma^s_t$ and $\widehat{\Gamma}^f_t$, $\widehat{\Gamma}^s_t$
in the current configuration and reference configuration, respectively.
The motion of the solid is followed in a Lagrangian way, therefore the domain $ \widehat{\Omega}^s $ is a Lagrangian domain and it is initially occupied by the solid we observe.
The domain $ \widehat{\Omega}^f $ is called ALE domain and is the domain on which we initially observe the fluid motion in a Eulerian way.
As a consequence of the solid movement, the domain on which we observe the fluid motion 
changes in time as well, so that we need to define a deformation for the fluid domain.
The domain $ \Omega^f_t $ is occupied only by fluid at each time $ t $. 
The moving domains $ {\Omega^f_t} $ and $ {\Omega^s_t} $ are called Eulerian domains. 
In order to describe the motion of the fluid and solid domains, let us define a $t$-parametrized family of
invertible and sufficiently smooth mappings $\mathcal{X}_t$ of the reference configuration $\widehat{\Omega}$ to the deformed
ones $\Omega_t$, so that 
\begin{align}
 \mathcal{X}_t & : \widehat{\Omega} \rightarrow \Omega_t \,,  \quad \mathcal{X}_t(\widehat{\mathbf{x}})  :=
 \widehat{\mathbf{x}} + {\bm{d}}(\widehat{\mathbf{x}},t)  \,.
\end{align}
The field ${\bm{d}}(\widehat{\mathbf{x}},t)$ is called displacement field.
The displacement field ${\bm{d}}(\widehat{\mathbf{x}},t)$ is determined separately in the fluid and solid parts 
as a solution of two different subproblems.
Its restrictions ${\bm{d}}^f(\widehat{\mathbf{x}},t)$ and ${\bm{d}}^s(\widehat{\mathbf{x}},t)$ 
are referred to as fluid domain displacement (or ALE displacement)
and solid displacement, respectively.
They are required to take on common values at the interface, namely
\begin{align}
  {\bm{d}}^s(\widehat{\mathbf{x}},t) = {\bm{d}}^f(\widehat{\mathbf{x}},t)   \,,\quad  \widehat{\mathbf{x}} \in \widehat{\Gamma}^i
  \,.
\end{align}
For every $\widehat{\mathbf{x}} \in \widehat{\Omega}$ and $t \ge 0$, %$ (\widehat{\mathbf{x}},t) \in \, Q_{\widehat{\Omega}} $, 
we also define
\begin{align}
 {\mathbf{F}} ({\bm{d}}(\widehat{\mathbf{x}},t) ) & = \widehat{\nabla}\mathcal{X}_t(\widehat{\mathbf{x}}) =
 I + \widehat{\nabla} {\bm{d}}(\widehat{\mathbf{x}},t) \,, \\ 
 J      ( {\bm{d}}(\widehat{\mathbf{x}},t) ) & = \det {\mathbf{F}}( {\bm{d}}(\widehat{\mathbf{x}},t) ) \,,\\
 {\mathbf{B}} ( {\bm{d}}(\widehat{\mathbf{x}},t) ) & = {\mathbf{F}} ( {\bm{d}}(\widehat{\mathbf{x}},t) )
 {\mathbf{F}}^T ( {\bm{d}}(\widehat{\mathbf{x}},t) ) \,,
\end{align}
where $\widehat{\nabla}$ refers to the gradient operator in the reference configuration.
The symbols  $ {\mathbf{F}} $ and  $ {\mathbf{B}} $ denote 
the deformation gradient tensor
and the left Cauchy-Green deformation tensor, respectively. In the following, 
$\nabla$ will refer to the gradient operator in the current configuration and the symbols $\widehat{\mathbf{n}}$ and ${\mathbf{n}}$ will denote the outward unit normal fields in the reference and in the current configuration, respectively.

\subsection{Weak Monolithic Formulation}\label{weak}
Here, the weak monolithic formulation of the FSI problem is described.
For the sake of simplicity, we denote with the same symbol $ ( \cdot, \cdot )  $ 
the standard inner products 
 either on $ L^2(\mathcal{O}) $, $ L^2(\mathcal{O})^n $
or $ L^2(\mathcal{O})^{n \times n} $,
for any open set $ \mathcal{O} \in \mathbb{R}^n$.
On the fluid-only boundary $ {\Gamma}^f_t $, 
we denote with $ {\Gamma}^f_{t,D} $ and $ {\Gamma}^f_{t,D,\bm{d}^f} $ the subsets of $ {\Gamma}^f_t $
on which Dirichlet boundary conditions on the velocity $\bm{u}^f $ 
and on the fluid domain displacement $\bm{d}^f $ are enforced, respectively.
Similarly, given the solid-only boundary $ {\Gamma}^s_t $, we denote with $ {\Gamma}^s_{t,D} $ the part
on which Dirichlet boundary conditions on the displacement $\bm{d}^s $ are enforced. 
In order to keep the exposition simple, we do not discuss the case of mixed Dirichlet-Neumann conditions
and the related definitions of function spaces and variational equations.
For any boundary subset $ \Gamma \subseteq  \partial\mathcal{O} $, we also denote with $ {H}^1_0(\mathcal{O};\Gamma) $
the subspace of functions in $ {H}^1(\mathcal{O}) $ with zero trace on $ \Gamma $.
Now define 
 \begin{align}
 & \bm{V}^{\bm{f}}  :=  \bm{H}^1(\Omega^f_t) \,, \quad
   \bm{V}^{\bm{s}}  :=  \bm{H}^1(\Omega^s_t)
   \,, \\ 
 & \bm{V}^{\bm{f}}_0              :=  \bm{H}^1_{0}(\Omega^f_t;         {\Gamma}^f_{t,D} ) \;, 
   \bm{V}^{\bm{s}}_0              :=  \bm{H}^1_{0}(\Omega^s_t;         {\Gamma}^s_{t,D} ) \;, 
   {\bm{V}}^{\bm{f}}_{0,\bm{d}^f} :=  \bm{H}^1_{0}({\Omega}^f_t; {\Gamma}^f_{t,D,\bm{d}^f} )
%    \widehat{\bm{V}}^{\bm{f}}_{0,\bm{d}^f} :=  \bm{H}^1_{0}(\widehat{\Omega}^f; \widehat{\Gamma}^f_{D,\bm{d}^f} )
   \,, \\ 
 & \bm{V} := 
  \{ \bm{v} = (\bm{v}^f,\bm{v}^s) \in 
  \bm{V}^{\bm{f}} \times \bm{V}^{\bm{s}} \text{ s. t. } \bm{v}^f = \bm{v}^s \text{ on } \Gamma^i_t \}  \,, \\
 & \bm{V}_0  := 
  \{ \bm{v} = (\bm{v}^f,\bm{v}^s) \in  
  \bm{V}^{\bm{f}}_0 \times \bm{V}^{\bm{s}}_0 \text{ s. t. } \bm{v}^f = \bm{v}^s \text{ on } \Gamma^i_t \} \,. 
 \end{align}
The mapping of $ {\bm{V}}^{\bm{f}}_{0,\bm{d}^f} $ to the reference domain 
is denoted as $\widehat{\bm{V}}^{\bm{f}}_{0,\bm{d}^f}$.

The weak monolithic FSI problem consists in finding 
$ (\bm{d}, \bm{u}, p) $ in 
$ \bm{V} \times \bm{V} \times L^2(\Omega_t) $ 
solution of a system that can be split into three parts:\\
the \textit{weak momentum balance}
\begin{align}
% solid momentum
 & \left( {\rho}^s  \frac{\partial {\bm{u}}}{\partial t} , {\bm{\phi}}^{m}       \right)_{{\Omega}_t^s}
 + \left(  {\bm{\sigma}}^s({\bm{d}},{p^s}) , {\nabla} {\bm{\phi}}^{m}  \right)_{{\Omega}_t^s} 
 - \left( {\rho}^s  {\bm{f}}^s    ,  {\bm{\phi}}^{m}      \right)_{{\Omega}_t^s}   
 %= 0\qquad \forall \; \bm{\phi}^{m} \in \bm{V}_0 , \label{momentum_sld} 
\nonumber \\
%fluid momentum
& + \left( \rho^f \frac{\partial \bm{u}}{\partial t} , \bm{\phi}^{m} \right)_{\Omega^f_t } 
  +  \left( \rho^f [ (\bm{u} - \frac{\partial \bm{d}}{\partial t} ) \cdot \nabla ] \bm{u} , \bm{\phi}^{m} \right)_{\Omega^f_t}
  \nonumber \\
& +  \left( \bm{\sigma}^f (\bm{u},p^f)     , \nabla \bm{\phi}^{m} \right)_{\Omega^f_t}
  -  \left( \rho^f \bm{f}^f   , \bm{\phi}^{m} \right)_{\Omega^f_t}
  = 0 
 \qquad
 \forall \; \bm{\phi}^{m} \in \bm{V}_0 , \label{momentum}
\end{align}
the \textit{weak mass continuity}
\begin{align}
% solid mass
\left( \widehat{J}( {\bm{d}} ) - 1 , {\phi}^{ps} \right)_{\widehat{\Omega}^s} 
 & = 0 
 \qquad
\forall \; {\phi}^{ps} \in L^2(\widehat{\Omega}^s),
\label{cntsld} \\ 
% fluid mass
 \left( \nabla \cdot  \bm{u}  , \phi^{pf} \right)_{{\Omega}^f_t} 
 & = 0
\qquad
\forall \; \phi^{pf} \in L^2({\Omega}^f_t) \,,
\label{cntfld}
\end{align}
and the \textit{weak kinematic equations}
 \begin{align} 
 % solid kinematic
 \left( {\bm{u}} - \frac{\partial {\bm{d}}}{\partial t}  , {\bm{\phi}}^{ks} \right)_{\widehat{\Omega}^s} 
  & = 0
 \qquad
 \forall \; {\bm{\phi}}^{ks} \in \bm{H}^1(\widehat{\Omega}^{s}), 
 \label{kinsld} \\  
 % fluid kinematic
  \left( k(\aled{\mathbf{x}})  \aled{\nabla} {\bm{d}}^f,  \aled{\nabla} {\bm{\phi}}^{kf}  \right)_{\widehat{\Omega}^f}
  & = 0 \qquad
 \forall \; {\bm{\phi}}^{kf} \in \bm{H}_0^1(\widehat{\Omega}^{f}; \widehat{\Gamma}^i) \cap \widehat{\bm{V}}^{\bm{f}}_{0,\bm{d}^f} \,. \label{kinfld}
 \end{align}
This formulation is known as a non-conservative ALE formulation. 
For a discussion about the conservative and non-conservative ALE formulation for a model problem of a scalar advection diffusion equation refer to \cite{formaggia2004stability}.
In the above formulation, equation \eqref{momentum} describes, in a monolithic form, the solid and fluid momenta, which are also referred to as the incompressible non-linear elasticity equation and the Navier-Stokes equation, respectively.
The symbols $\rho^f$ and $\rho^s$ denote the mass densities for the fluid and solid part, respectively, while
${{\mathbf{f}}}^f$  and ${{\mathbf{f}}}^s$ indicate the body force densities. 
The interface physical condition of normal stress continuity 
is enforced in the weak momentum balance, where the boundary integrals disappear due to the condition
\begin{equation}
   \bm{\sigma}^s (\bm{d},p^s) \mathbf{n}^s + \bm{\sigma}^f (\bm{u},p^f) \mathbf{n}^f = 0 \quad \text{ on } {\Gamma}^i_t\,. \label{interface}
\end{equation}
For the solid stress tensor  $ {\bm{\sigma}}^s $, we consider incompressible Mooney-Rivlin, 
whose Lagrangian description is given for every $\widehat{\mathbf{x}} \in \widehat{\Omega}^s$ and $t\ge 0$ by
\begin{align}
{\bm{\sigma}}^s ({\bm{d}} , p^s )
   & = - p^s  \mathbf{I} + 2 C_1 {\mathbf{B}}({\bm{d}} ) -  2 C_2 ( {\mathbf{B}}
   ({\bm{d}} ) )^{-1} \;, \label{Mooney-Rivelin-cauchy}
\end{align}
where the constants $C_1$ and $C_2$ depend on the mechanical properties of the material.
For the fluid stress tensor ${\bm{\sigma}}^f$, an incompressible Newtonian fluid is considered whose expression is given 
for every ${\mathbf{x}} \in \Omega^f_t$ and $t\ge 0$ by
\begin{align}  \label{newtonian}
{\bm{\sigma}}^f ( {\bm{u}} , p^f ) = - p^f {\mathbf{I}} + \mu ( \nabla {\bm{u}} + (\nabla {\bm{u}})^T ) \,,
\end{align}
where $\mu$ is the fluid viscosity.
Concerning velocity continuity, notice that the solid kinematic equation is just a change of variables 
without associated boundary conditions and its test functions are in $ \bm{H}^1(\widehat{\Omega}^{s}) $.
The interface velocity computed from this equation is an input for the fluid momentum balance.
With reference to displacement continuity,
notice also that the test functions in the weak fluid kinematic equation are in 
$ \bm{H}_0^1(\widehat{\Omega}^{f}; \widehat{\Gamma}^i) \cap \widehat{\bm{V}}^{\bm{f}}_{0,\bm{d}^f} $,
so that they vanish on the solid-fluid interface.
Thus, this equation does not affect the value of the 
displacement on the interface, which is an unknown of the problem 
that is evaluated by solving the other parts of the system.
% We remark that in deriving the weak FSI system we used the identity
% \begin{equation}
%  \left( \rho^f \frac{\partial \bm{u}}{\partial t}\biggl|_{\mathcal{\aleop}^f} , \bm{\phi} \right)_{\Omega^f_t } =
%   \frac{d}{dt} \left( \rho^f \bm{u} , \bm{\phi} \right)_{\Omega^f_t } 
%  -  \left( \rho^f  \left( \nabla \cdot \frac{\partial \bm{d}}{\partial t} \right)  \bm{u} , \bm{\phi} \right)_{\Omega^f_t} 
% \end{equation}
% that is proved in \cite{fernandez2005newton}.
% The advantage of this formula is that one need not compute time derivatives along the ALE deformations.

One of the novelties introduced by our ALE method for large structural displacements is the choice of the function $k(\aled{\mathbf{{x}}})$
in equation \eqref{kinfld}. The usual choice is a piecewise-constant function discontinuous across the element boundary, so that 
smaller elements in the mesh can be made stiffer. Here, the function $k(\aled{\mathbf{{x}}})$ is a distance function that assures an homogenous deformation throughout the entire mesh, not only for those elements close the valve leaflets. A more detailed discussion of the function $k(\aled{\mathbf{{x}}})$ can be found in section \ref{functionK}.

\section{Streamline Upwind Petrov Galerkin Stabilization}\label{supg}
This section describes the Streamline Upwind Petrov Galerkin (SUPG) method we use to stabilize our FSI simulations. 
The SUPG technique is very well known in the literature \cite{hughes1979, brooks1982streamline, le1993supg, tezduyar2006stabilization} and has been widely used for the Navier-Stokes equations \cite{franca1992stabilized, franca1993element, choi1997fractional, kirk2008development}.  
For each refinement level $k$, the stability parameter $\tau_k$ often depends on two constants: the mesh parameter $h_k$ and the inverse estimate constant $C_k$ \cite{franca1992stabilized}. An alternative definition of $\tau$ that obviates the use of $h_k$ and $C_k$ was proposed by Franca and Madureira in \cite{franca1993element}, where the stability parameters depend on the computation of the largest eigenvalue of a generalized eigenvalue problem.
In our SUPG stabilization technique, the stability parameter $\tau$ is determined as in \cite{franca1993element}.
Following our formulation of the problem, the SUPG method applied to the weak momentum balance \eqref{momentum} can be written as:
% $$\varphi_{GLS} = (u - v) \nabla \varphi\; \tau$$
% \\
find $ (\bm{d}, \bm{u}, p) $ in 
$ \bm{V} \times \bm{V} \times L^2(\Omega_t) $ satisfying\\
\newline
% \begin{equation} B_{SUPG}(u , p ; v ,\psi) = L_{SUPG}(v , \psi),\;\;\;\;
% \forall v \in V_h,\;\; \forall \psi \in Q_h, \end{equation}
% where 
% \begin{align}
% B_{SUPG}(u , p ; v ,\psi) & = \int_{\Omega} \dfrac{\partial u}{\partial t} v \, d\Omega + \int_{\Omega}\big[\big(u - \dfrac{\partial d}{\partial t}\big)\cdot \nabla \big]u v\, d\Omega+
% \mu \int_{\Omega} \big(\nabla u + (\nabla u)^T\big)\nabla v \, d\Omega - \int_{\Omega} p \nabla \cdot  v \, d\Omega+ \\ \nonumber
% & \int_{\Omega} \psi \nabla \cdot u \, d\Omega
% + \sum_{k=1}^{N} \int_{T_k} \big[ \dfrac{\partial u}{\partial t} + \big[\big(u - \dfrac{\partial d}{\partial t}\big)
% \cdot \nabla \big]u - \nabla \cdot \big( -pI + \mu \big(\nabla u + (\nabla u)^T\big) \big) \big] \cdot \\ \nonumber
% & \tau \big[ \dfrac{\partial  v}{\partial t} + \big[\big(u - \dfrac{\partial d}{\partial t}\big)
% \cdot \nabla \big] v - \nabla \cdot \big( -\psi I \big) \big] \, d\Omega 
% %+ \mu \big(\nabla  v + (\nabla  v)^T\big) \big) \big] \, d\Omega 
% \end{align}
%the \textit{weak momentum balance}
\begin{align}
% solid momentum
 & \left( {\rho}^s  \frac{\partial {\bm{u}}}{\partial t} , {\bm{\phi}}^{m}       \right)_{{\Omega}_t^s}
 + \left(  {\bm{\sigma}}^s({\bm{d}},{p^s}) , {\nabla} {\bm{\phi}}^{m}  \right)_{{\Omega}_t^s} 
 - \left( {\rho}^s  {\bm{f}}^s    ,  {\bm{\phi}}^{m}      \right)_{{\Omega}_t^s} \nonumber \\
%fluid momentum
& + \left( \rho^f \frac{\partial \bm{u}}{\partial t} , \bm{\phi}^{m} \right)_{\Omega^f_t } 
  +  \left( \rho^f [ (\bm{u} - \frac{\partial \bm{d}}{\partial t} ) \cdot \nabla ] \bm{u} , \bm{\phi}^{m} \right)_{\Omega^f_t}
  \nonumber \\
& + \left( \bm{\sigma}^f (\bm{u},p^f)      , \nabla \bm{\phi}^{m} \right)_{\Omega^f_t}
  - \left( \rho^f \bm{f}^f       , \bm{\phi}^{m} \right)_{\Omega^f_t} \nonumber \\
%SUPG
& + \sum_{k=1}^{N^f} \left( \rho^f \frac{\partial \bm{u}}{\partial t} , \bm{\phi}_{SUPG}^{m} \right)_{T_t^k} 
  + \sum_{k=1}^{N^f}  \left( \rho^f [ (\bm{u} - \frac{\partial \bm{d}}{\partial t} ) \cdot \nabla ] \bm{u} , \bm{\phi}_{SUPG}^{m} \right)_{T_t^k} \nonumber \\
& - \sum_{k=1}^{N^f}  \left( \nabla \cdot \bm{\sigma}^f (\bm{u},p^f) , \bm{\phi}_{SUPG}^{m} \right)_{T_t^k}
  - \sum_{k=1}^{N^f}  \left( \rho^f \bm{f}^f  , \bm{\phi}_{SUPG}^{m} \right)_{T_t^k}
  = 0 \qquad
 \forall \; \bm{\phi}^{m} \in \bm{V}_0. \label{momentum_supg}
\end{align}
% and \begin{equation}
% L_{SUPG}(v , \psi)=  \int_{\Omega} f \cdot  v \, d\Omega + \sum_{k=1}^{N} \int_{T_k} f \cdot  \tau 
% \big[ \dfrac{\partial  v}{\partial t} + \big[\big(u - \dfrac{\partial d}{\partial t}\big)
% \cdot \nabla \big] v - \nabla \cdot \big( -\psi I \big) \big] \, d\Omega 
% %+ \mu \big(\nabla  v + (\nabla  v)^T\big) \big) \big] \, d\Omega
% \end{equation}
The only difference between equations \eqref{momentum} and \eqref{momentum_supg} is in the stabilization terms present in the weak fluid momentum.
The set $\{T_{t}^k\}_{k=1}^{N^f}$ is a partition of $\overline{\Omega}_t^f$ such that $T_{t}^k \cap T_{t}^j = 0$ for $k\ne j$ and 
$\cup_{k=1}^{N^f} T_{t}^k = \overline{\Omega}^f_t$. Once a finite element discretization is introduced, 
$\{T_{t}^k\}_{k=1}^{N^f}$ will represent the finite element triangulation used, and $N^f$ will indicate the total number of fluid elements in the mesh. 
In \eqref{momentum_supg}, the function $\bm{\phi}_{SUPG}^m$ is defined as
\begin{equation}
\bm{\phi}_{SUPG}^m=\tau \left( \rho^f [ (\bm{u} - \frac{\partial \bm{d}}{\partial t} ) \cdot \nabla ] \bm{\phi}^m \right)\,,
\end{equation}
where the description of the parameter $\tau$ is given below.
Notice that the consistency between the fluid and solid region is automatically enforced with this stabilization technique. 
By definition of $\bm{\phi}_{SUPG}$, we have that $\bm{\phi}_{SUPG}=0$ in the solid, since on $\Omega^s_t$  
$\bm u-\frac{\partial \bm d}{\partial t}=0$. Namely, the solid velocity and the mesh velocity (or the time 
derivative of the solid displacement) are equal to each other. This is why in the solid momentum equation we omitted the contribution
from $\bm{\phi}_{SUPG}$.

For the stability parameter $\tau$, we consider the following design:
\begin{equation}\label{tausupg}
 \tau=\dfrac{\xi(\mbox{Re}_k(\mathbf{x}))}{\sqrt{\lambda_k}|\bm{u}(\mathbf{x})|_2\nu(\mathbf{x})}\;,
\end{equation}
where $\nu(\mathbf{x})$ indicates the Reynolds number and
\begin{align}
 & \xi(\mbox{Re}_k(\mathbf{x}))=
\begin{cases}
  \mbox{Re}_k(\mathbf{x})\,,\;  \text{ if }\;\;\;  0\leq \mbox{Re}_k(\mathbf{x})<1\\
  1\,,\;\;\;\;\;\;\;\;\;\;  \text{ if }\;\;\;  \mbox{Re}_k(\mathbf{x})\geq 1
\end{cases}\;,
\end{align}
\begin{equation}
 \mbox{Re}_k(\mathbf{x})=\dfrac{|u(\mathbf{x})|_2}{4\sqrt{\lambda_k}\nu(\mathbf{x})}\;,
\end{equation}
\begin{equation}
 \lambda_k= \mbox{ max}_{\bm{u} \in \bm{V}_0(T^k_t)} \dfrac{||\nabla \cdot (\nabla \bm{u} + (\nabla \bm{u})^T)||^2_{0,T^k_t}}{||(\nabla \bm{u} + (\nabla \bm{u})^T)||^2_{0,T^k_t}}\;,
\end{equation}
\begin{align}
 & |\bm{u}(\mathbf{x})|_2= 
% \begin{cases}\label{speed}
%   \sum_{i=1}^{n}\big(|u_i(\mathbf{x})|^p\big)^{1/p}\,,\;\;\text{ if }\;\;\; 1\leq p < \infty\\
%   \max\limits_{i=1,\cdots,n}|u_i(\mathbf{x})|\,,\;\;\text{ if }\;\;\; p=\infty 
% \end{cases}
  \sum_{i=1}^{n}\big(|u_i(\mathbf{x})|^2\big)^{1/2}\,. \label{speed}
%\end{cases}
%\;.
\end{align}
The parameter $\lambda_k$ is computed as the largest eigenvalue of the following generalized eigenvalue problem:
for $k=1, \dots, N^f$ find $\bm{w} \in \bm{V}_0(T_t^k)$ and $\lambda_k$ such that 
%$w \in (R_k(T_k)/\mathbb{R})$
\begin{equation}
(\nabla \cdot (\nabla \bm{w} +(\nabla \bm{w})^T), \nabla \cdot (\nabla \bm{\phi} +(\nabla \bm{\phi})^T))_{T_t^k} =
\lambda_k (\nabla \bm{w}, \nabla \bm{\phi})_{T_t^k},\;\;\forall \bm{\phi} \in \bm{V}_0(T^k_t).
\end{equation}
This problem is solved for the largest eigenvalue by the power method.\\
We remark that the present design of the stability parameters \eqref{tausupg}-\eqref{speed} does not require
explicit computations of inverse estimate constants, nor the computation of mesh parameters.
For every timestep $t_n$, the generalized eigenvalue problem is solved $s$ times, where $s$ indicates the number of mesh refinements.  
The parameters $\lambda_k$ are computed before the monolithic GMRES iteration. See section \ref{solver} for more details about the solver.

\section{Discretization} 
This section describes the discretization scheme adopted in our ALE approach to deal with large displacements of the valve leaflets. 
The proposed approach relies on the introduction of
\begin{itemize}
 \item a staggered in time mesh velocity in the discretization scheme to improve computational stability, 
 \item a scaling factor that measures the distance of a fluid element
 from the valve leaflets, to guarantee that there are no mesh entanglements in the fluid domain, and
 \item fictitious springs to model the contact force between closing valve leaflets.
\end{itemize}
In this section these three topics are discussed in details.

Recently, an extended ALE method to deal with large structural displacements based on a variational mesh optimization
technique has been proposed in \cite{basting2017extended}. In the aforementioned paper, a procedure for the alignment of
the structure interface with edges of the resulting triangulation is described. 
A method with similar properties was introduced in \cite{bejanov2008grid}, where the alignment procedure is based 
on explicit combinatorial considerations to approximate the interface using the fluid mesh. 

\subsection{Time Discretization Scheme with Staggered In Time Velocity}\label{staggered}
The weak FSI problem discussed in section \ref{supg} is now discretized in time and space.
% We construct the scheme by using backward finite differences for the time derivatives in the integrands,
% while we use a Crank-Nicolson scheme for the fluid and solid momentum balances. 
Let $ 0 = t_0 < t_1 < \dots < t_N = T $ be a subdivision of the time interval 
with constant time step $ \Delta t = t_{n+1} - t_{n}$. 
Let $ \Omega^f_{n} $ and $ \Omega^s_{n} $ be the fluid and solid domains at time $ t_n $. 
Let $(\bm{d}_n, \bm{u}_n, p_n)$ be the approximation at time $t_n$. 
For $ n = 0,1,\dots,N $, and for all test functions $\bm{\phi}^m \in \bm{V}_0$,
$\phi^{pf} \in L^2({\Omega}^f_t)$, $\phi^{ps} \in L^2(\widehat{\Omega}^s)$, 
${\bm{\phi}}^{{kf}} \in \bm{H}^1(\widehat{\Omega}^{s})$ and
$\bm{\phi}^{{ks}} \in \bm{H}_0^1(\widehat{\Omega}^{f}; \widehat{\Gamma}^i) \cap \widehat{\bm{V}}^{\bm{f}}_{0,\bm{d}^f}$
the time-discretized weak FSI problem consists in
finding the iterates $(\bm{d}_{n+1}, \bm{u}_{n+1}, p_{n+1}) $ for $ n = 0, 1, ..., N $ such that 

\textit{weak momentum balance}
\begin{align}\label{discrete_momentum}
% solid momentum
\sum_{ig}\Bigg(& \Big( \rho^s \frac{\left({\bm{u}}_{n+1} - {\bm{u}}_{n}\right)   }{\Delta t} , {\bm{\phi}}^{m} \Big)_{{\Omega}_{t_{ig}}^s}+
 \left( {\bm{\sigma}}^s({\bm{d}}(t_{ig}),p^s_{n+\frac{1}{2}}) , {\nabla} {\bm{\phi}}^{m}  \right)_{{\Omega}_{t_{ig}}^s} 
-  \Big({\rho}^s {\bm{f}}^s(t_{ig}),{\bm{\phi}}^{m} \Big) _{{\Omega}_{t_{ig}}^s} 
 \nonumber \\     
 & +\Big( \rho^f \frac{\left({\bm{u}}_{n+1} - {\bm{u}}_{n}\right)   }{\Delta t}, {\bm{\phi}}^{m} \Big)_{{\Omega}_{t_{ig}}^s}+
 \left( \rho^f \big[ (\bm{u}(t_{ig}) - \bm{v}(t_{ig})) \cdot \nabla \big] \bm{u}(t_{ig}) , \bm{\phi}^{m} \right)_{\Omega^f_{t_{ig}}} \nonumber \\  
 & +\left( {\bm{\sigma}}^f({\bm{u}}(t_{ig}),p^f_{n+\frac{1}{2}}) , {\nabla} {\bm{\phi}}^{m}  \right)_{{\Omega}_{t_{ig}}^f} 
   -  \Big({\rho}^f {\bm{f}}^f(t_{ig}),{\bm{\phi}}^{m} \Big) _{{\Omega}_{t_{ig}}^f} 
\nonumber \\
%SUPG
& + \sum_{k=1}^N \left( \rho^f \frac{(\bm{u}_{n+1}-\bm{u}_n)}{\Delta t} , \bm{\phi}_{SUPG}^{m} \right)_{T_{t_{ig}}^k} 
  + \sum_{k=1}^N  \left( \rho^f [ (\bm{u}(t_{ig}) - \bm{v}(t_{ig})) \cdot \nabla ] \bm{u}(t_{ig}) , \bm{\phi}_{SUPG}^{m} \right)_{T_{t_{ig}}^k} \nonumber \\
& - \sum_{k=1}^N  \left( \nabla \cdot \bm{\sigma}^f (\bm{u}(t_{ig}),p_{n+\frac{1}{2}}^f) , \bm{\phi}_{SUPG}^{m} \right)_{T_{t_{ig}}^k}
  - \sum_{k=1}^N  \left( \rho^f \bm{f}^f(t_{ig})  , \bm{\phi}_{SUPG}^{m} \right)_{T_{t_{ig}}^k}
\Bigg) w_{ig}  \Delta t
=0\,,
\end{align}

\textit{weak mass continuity}
\begin{equation}\label{cont}
% solid mass
\left( \widehat{J}( {\bm{d}}_{n+1} ) - 1 , {\phi}^{ps}  \right)_{\widehat{\Omega}^s} 
% fluid mass
+ \left( \nabla \cdot  \bm{u}_{n+1}, \phi^{pf} \right)_{{\Omega}^f_{n+1}} = 0\,,
\end{equation} 

\textit{weak kinematic equation}
\begin{align}
% solid kinematic
\left( {\bm{u}}_{n+1} - \frac{{\bm{d}}_{n+1} -{\bm{d}}_{n}}{\Delta t} , {\bm{\phi}}^{{ks}} \right)_{\widehat{\Omega}^s} 
% fluid kinematic
+ \left( k(\aled{\mathbf{x}})  \aled{\nabla} {\bm{d}}_{n+1},
\aled{\nabla} {\bm{\phi}}^{{kf}}  \right)_{\widehat{\Omega}^f} = 0\,. \label{discrete_kin}
\end{align}
For both the fluid and solid momentum balances we sum over the Gauss points in the interval from $t_n$ to $t_{n+1}$, where
the symbol $ig$ indicates the Gauss points. We identify such Gauss points in the interval $(t_n, t_{n+1})$ using a mapping from the  $(-1, 1)$ to $(t_n, t_{n+1})$; namely a Gauss point in $(t_n, t_{n+1})$ is given by
\begin{equation}
 t_{ig} = t_n + \dfrac{x_{ig}+1}{2} \Delta t \,\,,
\end{equation}
where $x_{ig}$ corresponds to a Gauss point in $(-1, 1)$. In equation \eqref{momentum_supg}, $\omega_{ig}$ is the weight corresponding to the Gauss point $x_{ig}$ in $(-1, 1)$.
All accelerations in the momentum balance are considered constant 
in the interval from $t_n$ to $t_{n+1}$ and evaluated  as 
\begin{equation}
 {\bm{a}(t)} = \dfrac{{\bm{u}}_{n+1} - {\bm{u}}_{n}}{\Delta t}\; \mbox{ for all } t \in (t_n, t_{n+1}).
\end{equation}
Similarly, on each element the pressure $p$ is considered constant in the interval $(t_n, t_{n+1})$. 
Namely $p$ is piecewise constant in time throughout the entire simulation for both the solid and the fluid part.

With $\bm{v}$ we indicate the mesh velocity, that is computed following a linear interpolation in the time interval $(t_n, t_{n+1})$, namely
\begin{equation}
 \bm{v}(t_{ig})=(1-s)\bm{v}_{n} + s\, \bm{v}_{n+1}\;,
\end{equation}
where $s=\frac{t_{ig}-t_n}{\Delta t}$.
The value $\bm{v}_{n+1}$ is computed using the staggered approach 
\begin{equation}
\bm{v}_{n+1} = \frac{\bm{d}_{n+1}-\bm{d}_n}{\Delta t} \quad \mbox{ for } n\ge 1\,.
\end{equation}
Also in the kinematic equation, a staggered in time velocity is considered in the solid region. Moreover, since the first term of \eqref{discrete_kin} does not involve any space derivative, the corresponding mass matrix arising from the discretization is replaced with a lumped mass matrix.
Overall, the staggered approach and the use of a lumped mass matrix guarantee more stability.

Note that the discrete continuity equation \eqref{cont} is always solved at the end of the interval, meaning at $t_{n+1}$.
% To compute the mesh velocity $\bm{v}^{n+1}$ we adopt a staggered approach. 
% \begin{equation}
%  \frac{1}{2}\left(\bm{v}^{n+1}+\bm{v}^{n}\right)=\dfrac{\bm{d}_{n+1}-\bm{d}_{n}}{\Delta t}\label{standardDerivative}
% \end{equation}
% but we approximate it with $$\bm{v}^{n+1}=\dfrac{\bm{d}_{n+1}-\bm{d}_{n}}{\Delta t}$$ meaning we consider $u_{mesh}^{n+1} \approx u_{mesh}^{n+1/2}$.
% This staggered approach is again used to compute the solid velocity, while to compute the fluid velocity we rely on \eqref{standardDerivative}.\\

\subsection{Geometry and Aligned Mesh}\label{functionK}
\begin{figure}
 \begin{center}
%   (a) \includegraphics[scale=0.9]{figures/valve_mesh_zoom} \hspace{1cm}
%   (b) \includegraphics[scale=0.9]{figures/valve_tip}
  (a) \includegraphics[scale=0.5]{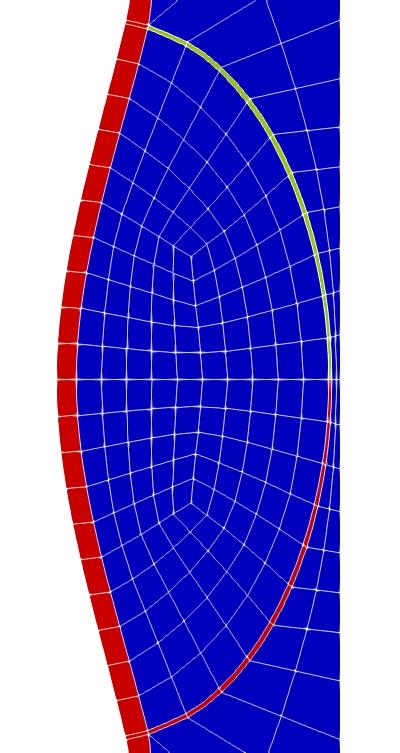} \hspace{0.5cm}
  (b) \includegraphics[scale=0.5]{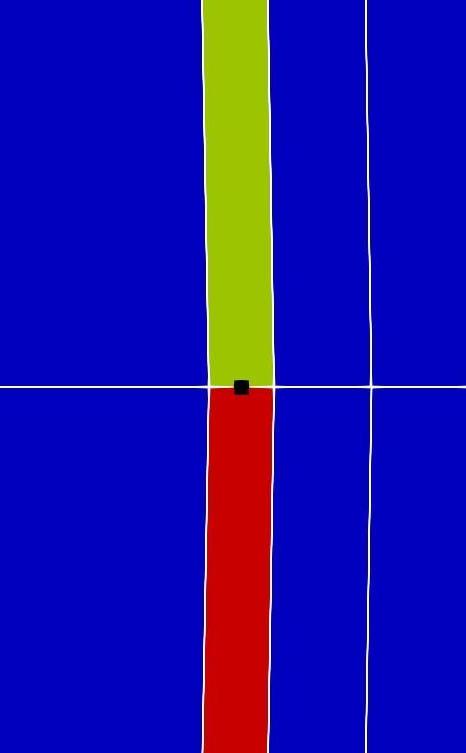}
  \caption{(a) zoom of the coarse mesh highlighting the \textit{fluid leaflet}, in green; 
  (b) zoom of the leaflet tip, where the black square indicates the location of the middle point $m$.}\label{zoom_n_tip}
 \end{center}
\end{figure}

The function $k(\aled{\mathbf{x}})$ in equation \eqref{discrete_kin} and \eqref{kinfld} is usually chosen to be a piecewise-constant function discontinuous across the element boundary so that smaller elements in the mesh can be made stiffer.
In this work, we make a different choice for the function $k$, a choice more suitable for the problem under consideration.
To better describe the function $k$, let's first focus on the geometry used to perform the venous valve simulations.
For our 2D geometry, half of a blood vessel is considered, namely only the motion of one of the two venous valve leaflets is analyzed. The other half of the geometry, together with its mesh, can be easily reconstructed by vertical symmetry. 
\begin{figure}[!b]
 \begin{center}
  (a) \includegraphics[scale=0.6]{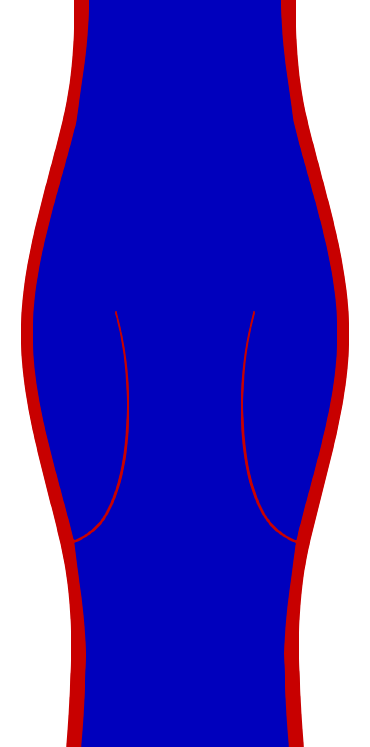} \hspace{1cm}
  (b) \includegraphics[scale=0.6]{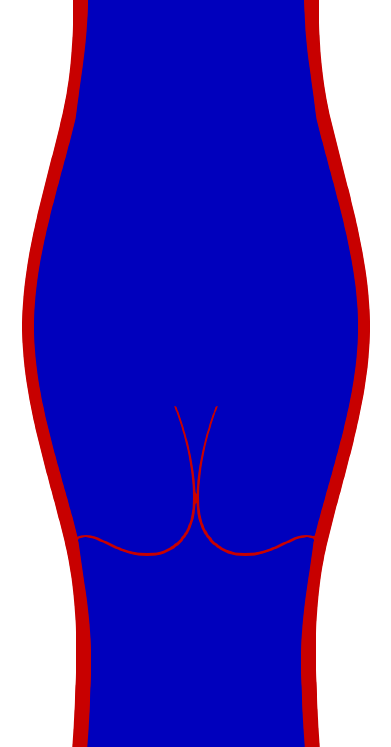} 
  \caption{Vein valve geometry, where the red color indicates the solid domain and the blue color indicates the fluid domain;
  (a) open valve configuration (b), closed valve configuration.}\label{entire_valve}
   \end{center}
\end{figure}
Figure \ref{entire_valve} (a) shows the entire geometry (reconstructed by vertical symmetry) in the equilibrium phase, where the red color indicates the solid domain and the blue color indicates the fluid domain. Figure \ref{entire_valve} (b) shows the entire geometry in the closed phase.
The upper part of the geometry in Figure \ref{entire_valve} (a), together with its mesh,
has been created by symmetry with respect to the $x$-axis; therefore, the leaflet geometry is present in the fluid domain as well.
This guarantees a more uniform mesh for the entire geometry, since the mesh in the upper part is exactly the same as the one in the lower part, but mirrored. The leaflet geometry present in the fluid domain will be called, from now on, \textit{fluid leaflet}. Figure \ref{zoom_n_tip} (a) shows a zoom of the fluid leaflet.
Given a fluid element in the mesh and being $\aled{\mathbf{x}}$ the center point of this element, the distance function $k(\aled{\mathbf{x}})$ is defined as
\begin{equation}
 k(\aled{\mathbf{x}}) = \dfrac{1}{1+c\cdot k_1(\aled{\mathbf{x}})}\label{scaleWeightNoJacobian}
\end{equation}
where $c \in \mathbb{R}^{+}$, $c\ge 1$ and $k_1(\aled{\mathbf{x}})$ is a function that measures the distance from $\aled{\mathbf{x}}$, the center of the fluid element, to $m$, the middle point at the tip of the valve leaflet. 
The location of $m$ in our 2D configuration is shown in Figure \ref{zoom_n_tip} (b). Thus, $k_1(\aled{\bm{{x}}})$ has the following expression
\begin{equation}
 k_1(\aled{\mathbf{x}}) = \sqrt{(\aled{\mathbf{x}}_1 - m_1)^2+(\aled{\mathbf{x}}_2 - m_2)^2}
\end{equation}
where the pair $(\aled{\mathbf{x}}_1, \aled{\mathbf{x}}_2)$ indicates the coordinates of $\aled{\mathbf{x}}$, and the pair $(m_1, m_2)$ indicates the coordinates of $m$.
The function $k(\aled{\mathbf{x}})$ guarantees that there are no mesh entanglements in the fluid domain by making the elements close the leaflet stiffer and the elements close to the vein wall softer. 

Another use of the function $k(\aled{\mathbf{x}})$ is in the modeling the movement of the leaflet meshed in the fluid domain. 
Our goal is that the \textit{fluid leaflet} follows the movement of the actual leaflet, so that there is uniformity of the mesh at every instant of time $t$, not only at $t=0$.
To impose this behavior, the \textit{fluid leaflet} has to be made stiffer, otherwise it would simply follow the fluid motion.
For this reason, for every element that composes the \textit{fluid leaflet}, the function $k(\aled{\mathbf{x}})$ defined in \eqref{scaleWeightNoJacobian} is multiplied by a constant $a \in \mathbb{R}^{+}$, $a\ge 1$. Namely, for every point ${\aled{\mathbf{x}}}_{FL}$, we define
\begin{equation}
 k({\aled{\mathbf{x}}}_{FL}) = \dfrac{a}{1+c\cdot k_1({\aled{\mathbf{x}}}_{FL})}\,.
\end{equation}
where ${\aled{\mathbf{x}}}_{FL}$ denotes the center of an element composing the \textit{fluid leaflet}.
In the numerical tests described in section \ref{tests}, we consider $a=100$ and $c=10,000$.\\

\subsection{Fictitious Springs}\label{springs}
As Figure \ref{entire_valve} shows, for our 2D geometry, half of a blood vessel is considered, namely only the motion of one of the two venous valve leaflets is analyzed. 
To make sure that the vein leaflet does not exit the domain (i.e does not overcome the axis of symmetry), we impose fictitious elastic springs on the fluid boundary throughout the simulation. With this technique, the leaflet never touches the boundary but can get very close to it. The closeness to the boundary is controlled by the given elastic parameters. The stiffness of each fictitious elastic spring is such that the spring carries a minimum force that allows the leaflet not to exit the fluid domain. 

Let $\Omega_t^M$ be the subregion of $\Omega^f_t$ belonging to the meniscus domain between the valve leaflet and the axis of symmetry. 
Each fictitious spring contributes with a force given by
\begin{equation}
 E(\bm{x}(t)) = E_0 \cdot \text{e}^{\frac{ x(t) - x_0 }{ h }}\; \text{ with } x(t) \le 0,
\end{equation}
where $x$ represents is the first component of $\mathbf{x}$.
The elasticity constant $E$ contributes only for $x(t)$ close enough to $x_0$
and it decays exponentially fast for $x<x_0$.
Notice that the non-linearity of $E$ comes both from the fact that $E$ is defined by an exponential function
and from its dependance on $\mathbf{x}:=\widehat{\mathbf{x}} + {\bm{d}}(\widehat{\mathbf{x}},t)$.
Considering the stabilized weak momentum balance illustrated in section \ref{supg},
the fictitious spring force is added in the momentum equation \eqref{momentum_supg} using the following 
\begin{equation}\label{spring_force}
 \Big( E(\mathbf{x}(t)) \nabla \bm{d}(t), \nabla \bm{\phi}^m \Big)_{\Omega^M_t}\,.
\end{equation}
The values considered in our numerical tests are: $E_0=$1.e-03, $x_0=-$1.0e-05 and $h$=1e-03.
Notice that the smaller the elastic constants become, the more the leaflet gets close to the boundary. 

This technique is widely used in structural dynamics and for the resolution of dynamics contact problems 
\cite{gambhir2013stability, filiatrault2013elements}. The time discretized version of \eqref{spring_force} that has to be added to equation \eqref{discrete_momentum} is given by
\begin{equation}
 \sum_{ig}\left( \sum_{k=1}^{S} 
  \Big( E(\mathbf{x}(t_{ig})) \nabla \bm{d}(t_{ig}), \nabla \bm{\phi}^m \Big)_{M_{t_{ig}}^k}
 \right)\omega_{ig} \Delta t\,,
\end{equation}
where $\{M_{t_{ig}}^k\}_{k=1}^{S}$ is a partition of $\Omega^M_{t_{ig}}$ such that
$M_{t_{ig}}^k \cap M_{t_{ig}}^j = 0$ for $k\ne j$ and 
$\cup_{k=1}^N M_{t_{ig}}^k = \overline{\Omega}^M_{t_{ig}}$.

\section{Monolithic Newton-Krylov solver}\label{solver}
In this section the Newton-Krylov solver employed in the solution of the monolithic FSI system is briefly described.
Once the weak monolithic formulation presented in section \ref{staggered} is discretized 
in space using appropriate finite element spaces,
the corresponding Jacobian matrix is obtained by an exact Newton linearization implemented by automatic differentiation 
(\cite{hoganadept}). The solution of the linear systems is performed using a GMRES solver
preconditioned by a geometric multigrid algorithm.  The smoother is of modified Richardson type, 
in turn preconditioned by a restricted additive Schwarz method.
The coarse grid correction problem is dealt with by a direct solver of the monolithic system.
Since the core of this work is the discretization scheme presented in the previous section, it is not our intention to provide a full account of the solver. For more details about the solver refer to \cite{FSIarticle} and \cite{COMPDYN2017}.

\subsection{Structure of the Jacobian} Let ${\bm{J}}^{(k)}$ denote the 
exact Jacobian at a non-linear step $k$. Ordering the variables as 
\begin{align}
\begin{bmatrix}
{\bm{d}}^{s}     \,
{\bm{d}}^{i}     \,
{\bm{d}}^{f}     \,
\vline         \,
{\bm{u}}^{s}     \,
{\bm{u}}^{i}     \,
{\bm{u}}^{f}     \,
\vline         \,
p^{s}          \,
p^{f}	       \,
\end{bmatrix}^\intercal
\end{align}
we consider the Jacobian to have the following block structure \begin{align}
{\bm{J}}^{(k)} =
  \begin{bmatrix}[ ccc || ccc || cc ] \label{jacobian}
 S_{\bm{d}^s}^{\bm{d}^s}  & S_{\bm{d}^i}^{\bm{d}^s} &  0  & S_{\bm{u}^s}^{\bm{d}^s} & S_{\bm{u}^i}^{\bm{d}^s} & 0 &  S_{p^s}^{\bm{d}^s}  &  0\\
 I_{\bm{d}^s}^{\bm{d}^i}  & I_{\bm{d}^i}^{\bm{d}^i} &  I_{\bm{d}^f}^{\bm{d}^i}  & I_{\bm{u}^s}^{\bm{d}^i} & 
 I_{\bm{u}^i}^{\bm{d}^i} & I_{\bm{u}^f}^{\bm{d}^i} &  I_{p^s}^{\bm{d}^i}  &  I_{p^f}^{\bm{d}^i}\\
 0  &  A_{\bm{d}^i}^{\bm{d}^f} & A_{\bm{d}^f}^{\bm{d}^f}& 0 & 0  &  0 &  0  &  0\\
\hline
\hline
K_{\bm{d}^s}^{\bm{u}^s} & K_{\bm{d}^i}^{\bm{u}^s}   &   0  &  K_{\bm{u}^s}^{\bm{u}^s}  &  K_{\bm{u}^i}^{\bm{u}^s}  & 0 &   0   &   0\\
K_{\bm{d}^s}^{\bm{u}^i} & K_{\bm{d}^i}^{\bm{u}^i}   &   0  &  K_{\bm{u}^s}^{\bm{u}^i}  &  K_{\bm{u}^i}^{\bm{u}^i}  & 0 &   0   &   0\\
0 & F_{\bm{d}^i}^{\bm{u}^f }  & F_{\bm{d}^f}^{\bm{u}^f }   &  0   &  F_{\bm{u}^i }^{\bm{u}^f } &F_{\bm{u}^f }^{\bm{u}^f } &   0  &
F_{p^f}^{\bm{u}^f }\\
\hline
\hline
V_{\bm{d}^s}^{p^s} &  V_{\bm{d}^i}^{p^s} & 0   &   0 & 0 & 0  &   0  &  0\\
0  & W_{\bm{d}^i}^{p^{f}} &W_{\bm{d}^f}^{p^{f}}  &   0  &  W_{\bm{u}^i }^{p^{f}} & W_{\bm{u}^f }^{p^{f}} &  0 & 0\\
 \end{bmatrix}
\begin{matrix}[l]
\vspace{2pt} \mbox{Momentum Solid} \\ 
\vspace{2pt} \mbox{Momentum Interface} \\
\vspace{6pt} \mbox{Kinematic fluid} \\ 
\vspace{2pt} \mbox{Kinematic Solid} \\
\vspace{2pt} \mbox{Kinematic Interface} \\
\vspace{6pt} \mbox{Momentum Fluid} \\ 
\vspace{2pt} \mbox{Continuity Solid} \\
\vspace{6pt} \mbox{Continuity Fluid}
\end{matrix}
\end{align}
where we use the symbols $K$ for the kinematic equation in the solid and on the 
fluid-solid interface \eqref{kinsld}, $A$ for the kinematic ALE displacement equation in the fluid \eqref{kinfld}, 
$S$ and $F$ for the momentum equation in the solid and fluid respectively \eqref{momentum},
%$S$ for the momentum equation in the solid \eqref{momentum_sld},
%$F$ for the momentum equation in the fluid \eqref{momentum_fld}, 
$I$ for the momentum equation \eqref{interface} on the fluid-solid interface,
$V$ for the continuity equation \eqref{cntsld} in the solid and 
$W$ for the continuity equation \eqref{cntfld} in the fluid.

The equations and unknowns are ordered following a field-ordering approach as in \cite{fernandez2005newton}.
It is important to notice that different orderings, though equivalent mathematically, can have a significant effect on the 
convergence properties and computational time of the solver, especially in the parallel setting.

\subsection{Geometric Multigrid Preconditioner}
As a preconditioner to the outer monolithic GMRES iteration, we consider a Geometric Multigrid algorithm.
Let $ {\bm{\Phi}}({\Omega}_{h_{l}}) $ and $ {\bm{\Psi}}({\Omega}_{h_{l}}) $
be the finite element spaces associated with each level of the triangulation $ \{\Omega_{h_l}\}_{l=1}^{L} $
with relative mesh sizes $ h_l $. The \textit{prolongation} $I^{l}_{l-1}$ and \textit{restriction} $I_{l}^{l-1}$ operators
\begin{align}
 I_{l-1}^{l} & :
 \bm{\Phi}({\Omega}_{h_{l-1}}) \times 
 \bm{\Phi}({\Omega}_{h_{l-1}}) \times 
 {\Psi}({\Omega}_{h_{l-1}})
 \rightarrow 
 \bm{\Phi}({\Omega}_{h_{l}}) \times 
 \bm{\Phi}({\Omega}_{h_{l}}) \times 
 {\Psi}({\Omega}_{h_{l}})   \,, 
 \\
 I_{l}^{l-1} & :
 \bm{\Phi}({\Omega}_{h_{l}}) \times 
 \bm{\Phi}({\Omega}_{h_{l}}) \times 
 {\Psi}({\Omega}_{h_{l}})    
 \rightarrow 
 \bm{\Phi}({\Omega}_{h_{l-1}}) \times 
 \bm{\Phi}({\Omega}_{h_{l-1}}) \times 
 {\Psi}({\Omega}_{h_{l-1}})
\end{align}
are defined, respectively, as the natural injection from the coarse to the fine space and the adjoint 
of $ {I}^{l}_{l-1} $ with respect to the $ L^2 $ inner product. Clearly, the matrix representations 
of these operators, $ {\bm{I}}_{l-1}^{l} $ and $ {\bm{I}}_{l}^{l-1} $, depend on the block row ordering
of the Jacobian \eqref{jacobian}, so that a different Jacobian structure affects the structure of 
$ {\bm{I}}_{l-1}^{l} $ and $ {\bm{I}}_{l}^{l-1} $.
The block structures of the prolongation and restriction operators are
\begin{equation*} 
 {\bm{I}}_{l}^{l-1} =
\begin{bmatrix}[ ccc || ccc || cc]
   R_{\bm{d}^s}^{\bm{d}^s}          &    0      &      0  &   0   &       0        & 0   &   0     &      0         \\
   R_{\bm{d}^{s}}^{\bm{d}^{i}}    &    R_{\bm{d}^{i}}^{\bm{d}^{i}}    & 0   &   0   &
   0  &   R_{\bm{u}^{f}}^{\bm{d}^{i}}    &    0   &   0\\
   0   &   0    &  R_{\bm{d}^f}^{\bm{d}^f}  &         0        &         0          &     0    &    0 &     0  \\
\hline
\hline
   0   &   0              &   0            &    R_{\bm{u}^s}^{\bm{u}^s}  &    0     &     0 
   &   0          &     0         \\
   0   &   0              &   0            &    R^{\bm{u}^{i}}_{\bm{u}^s}   &
   R_{\bm{u}^{i}}^{\bm{u}^{i} }  &    0  &   0 & 0 \\
   0   &  0   &   0    &         0    &    0  & 
   R_{\bm{u}^f }^{\bm{u}^f }  &   0          &     0      \\
\hline
\hline
   0            &    0             &         0          &       0  &    0
   &     0              &     R_{p^s}^{p^s} & 0    \\
   0          &   0  & 0   &  0 &   0          &     0  
   &     0      &       R_{p^{f}}^{p^{f}} \\
\end{bmatrix}
\,,\;\;
 {\bm{I}}_{l-1}^{l} =
\begin{bmatrix}[ ccc || ccc || cc]
   P_{\bm{d}^s}^{\bm{d}^s}   &    P_{\bm{d}^{s}}^{\bm{d}^{i}}      &      0  &   0   &       0        & 0   &   0     &      0         \\
   0    &    P_{\bm{d}^{i}}^{\bm{d}^{i}}    & 0   &   0   &   0  &   0    &    0   &   0\\
   0   &   P_{\bm{d}^{f}}^{\bm{d}^{i}}    &  P_{\bm{d}^f}^{\bm{d}^f}  &    0        &     0     &     0    &    0 &     0  \\
\hline
\hline
   0   &   0              &   0            &    P_{\bm{u}^s}^{\bm{u}^s}  &    P^{\bm{u}^{i}}_{\bm{u}^s}     &     0 
   &   0          &     0         \\
   0   &   0              &   0            &    0   &   P_{\bm{u}^{i}}^{\bm{u}^{i} }  &    0  &   0 & 0 \\
   0   &  0   &   0    &         0    &    P_{\bm{u}^f }^{\bm{u}^i }  & 
   P_{\bm{u}^f }^{\bm{u}^f }  &   0          &     0      \\
\hline
\hline
   0            &    0             &         0          &       0  &    0
   &     0              &     P_{p^s}^{p^s} & 0    \\
   0          &   0  & 0   &  0 &   0          &     0  
   &     0      &       P_{p^{f}}^{p^{f}} \\
\end{bmatrix}
\,.\\
\end{equation*}

\subsection{Richardson-Schwarz Smoother} In the smoothing process, we first partition the whole domain 
into the fluid and solid subregions, 
and then we further divide each subregion into smaller non-overlapping blocks ${\Omega}_k, k = 1, ...,N$.
On each subdomain ${\Omega}_k$, we construct a subdomain preconditioner $\bm{B}_k$, 
which is a restriction of the Jacobian matrix $\bm{J}$; that is, 
it contains entries from $\bm{J}$ corresponding to the degrees of freedom (DOFs) 
contained in the corresponding subdomain ${\Omega}_k$.
The exchange of information between blocks is guaranteed by the fact that the support of the
test function associated to the displacement and velocity DOFs extends to the neighboring elements.
The restricted version of the additive Schwarz (AS) preconditioner used in
the Richardson scheme for the FSI Jacobian system is
\begin{equation}
 {\bm{B}}^{-1} = \sum_{k=1}^N ({\bm{R}}^0_k)^T {\bm{B}}_k^{-1}({\bm{R}}^{\delta}_k)\,,
\end{equation}
where $\bm{R}_k$ indicates a restriction matrix that maps the global vector of degrees of freedom 
to those belonging to the subdomain $\Omega_k$, and 
$\bm{R}^{0}_k$ is a restriction matrix that does not include the overlap, while $\bm{R}^{\delta}_k$ does.

\section{Numerical Results}\label{tests}
\begin{figure}
\begin{center}
 \includegraphics[scale=1.0]{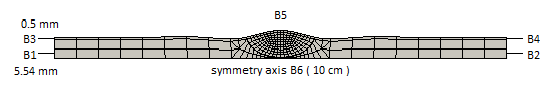}
 \caption{Geometry with coarse mesh, lengths and boundary names.}\label{valve_entire2}
\end{center}
\end{figure}
In this section, the numerical tests aimed at assessing our ALE approach are described. 
The solver and all the tests have been implemented in the in-house
finite element C++ library FEMuS (\url{https://github.com/FeMTTU/femus}). The solver has
been implemented using the GMRES solver and the geometric multigrid preconditioner interface
implemented in the PETSc toolkit \cite{balay2012petsc}.
For our simulations, a 2D geometry of a venous valve is used. As already mentioned in section \ref{functionK}, half of a blood vessel is considered, namely only the motion of one of the two venous valve leaflets is analyzed. Figure \ref{valve_entire2} shows the entire 2D configuration used for the tests, together with its mesh, lengths and the names given to the different boundaries.
In this geometry, we identify 6 different boundaries: $B1$ and $B2$ represent the fluid boundaries at the bottom and top of the vein, respectively, $B3$ and $B4$ represent the solid boundaries at the bottom and top of the vein, respectively,
$B5$ indicates the lateral solid boundary and $B6$ is the symmetry axis.
The vein has a lumen diameter of $5.54$ mm and a length of $10$ cm.
The thickness of the vein wall is $0.5$ mm. For the valve leaflet a thickness of $0.065$ mm is specified.
%$0.05$ mm is specified, as in \cite{narracott2010analysis}. 
All solid materials are modeled as incompressible non-linear elastic, as explained in section \ref{weak}.
Both the vein wall and the valve leaflet are considered to have the same density, $960$ kg/m$^3$, and Poisson's ratio, $0.5$.
Blood is considered a Newtonian fluid with a density of $1060$ kg/m$^3$ and a viscosity of $2.2 \times 10^{−3}$ Pa$\cdot$s. Although blood is known to be non-Newtonian in general, several studies, such as
\cite{bazilevs2010computational, crosetto2011fluid, turek2011numerical},
assume it to be Newtonian, as we do in this paper.

Below, two tests using the described 2D valve model are performed. Test $1$ studies the mesh and process independence of the proposed ALE method and solver, and Test $2$ investigates the sensitivity of valve dynamics to the elastic modulus of the valve leaflets. The wall and leaflet Young's moduli will be specified in each test case. For both tests, the same boundary conditions are used. At the bottom fluid boundary $B1$, shown in Figure \ref{valve_entire2}, a normal stress boundary condition of the form
\begin{align}
\boldsymbol (\bm{\sigma}(\bm{u},p^f) \cdot \mathbf n) \cdot \mathbf n &= 15 \sin(2\pi t ) \quad \mbox{[Pa]}, \nonumber \\
\bm u  \cdot \boldsymbol \tau &= 0, \nonumber \\
\bm{d} &= \bm{0}.
\end{align}
is specified, where $\bm{\tau}$ indicates the tangential vector to the boundary.
At the top boundary $B2$, the same condition has been applied but with the opposite sign. 
The vein is considered clamped, so at $B3$ and $B4$ the following condition is applied 
\begin{align}
 \bm{d} = \bm{0}.
\end{align}
At the lateral solid boundary $B5$, a stress boundary condition of the form
\begin{align}
 \bm{\sigma}(\bm{d},p^s) \cdot \mathbf n &= \bm{0}.
\end{align}
is specified. Finally, at the symmetry axis $B6$, we require 
\begin{align}
\bm u  \cdot \mathbf n &= 0, \nonumber \\
\bm d  \cdot \mathbf n &= 0, \nonumber \\
\dfrac{\partial \mathbf u}{\partial \boldsymbol\tau} &= 0, \nonumber \\
\dfrac{\partial \mathbf d}{\partial \boldsymbol\tau} &= 0.
\end{align}

\subsection{Test 1: Mesh and Process Independence}\label{test1}
In Test $1$, the mesh and process independence of the proposed ALE method and solver are studied.
A Young's modulus of $260$ MPa for the vein wall and of $1.5$ MPa for the valve leaflet are considered.
These parameters do not represent realistic biological values, which is understandable since we are performing 2D simulations that require much stiffer structure properties in order to have reasonable deformations.
To investigate mesh independence, we consider three uniform refinement levels starting from a coarse mesh that has $292$ total elements, namely we consider $3$, $4$ and $5$ levels of refinement. The three simulations have been performed in parallel using 36 processes. To study process independence, we consider the $4$ refinements level case and perform this simulation using $18$ and $72$ processes as well. 
All data are collected considering one single valve period, which we define as $1$ s. 
The time step considered is $\dfrac{1}{64}$s, meaning that one valve period is composed of 64 iterations.

To quantify the dependance of the solver on the mesh, for every non-linear step $s$ at time step $t_i$, $i=1,\dots,64$, 
let's define the average convergence rate in the linear solvers 
as $\rho_{i,s} = \left(\dfrac{r_n}{r_0}\right)^N_{i,s}$, where $N_{i,s}$ is the number of linear steps (in each non-linear step $s$) and $r_n$ and $r_0$ represent the residuals of two subsequent iterations. %at time $t_n$ and $t_0$.
Over a single valve period, let's define two quantities, $N$ and $\rho$, as 
\begin{align}
 N&=\dfrac{\sum_{i=1}^{64}\sum_{s=1}^{s_i^{max}} N_{i,s}}{\sum_{i=1}^{64}\sum_{s=1}^{s_i^{max}}1},\\
 \rho&=\dfrac{\sum_{i=1}^{64}\sum_{s=1}^{s_i^{max}} \rho_{i,s}}{\sum_{i=1}^{64}\sum_{s=1}^{s_i^{max}}1},
\end{align}
where $s_i^{max}$ is the maximum number of non-linear steps at time step $t_i$, $i=1,\dots,64$.
The symbol $\rho$ indicates the average of the $\rho_{i,s}$'s over a single valve period, while $N$ represents the average of the $N_{i,s}$'s over a single valve period. The values of $N$ and $\rho$ for the three refinement levels considered, namely 
$3$, $4$ and $5$ levels of refinement, are shown in Table \ref{mesh}.
\begin {table}[!h] %[!tpb]
\setlength\tabcolsep{5.5pt} % default value: 6pt
\begin{center}
	\begin{tabular}{|c|c|c|c|c|} \hline	
	\multirow{2}{*} {}   & \multicolumn{3}{|c|}{Mesh Independence of the Solver} \\ \cline {2-4}
	& \multicolumn{1}{|c|} {$s^{max}$ } & \multicolumn{1}{|c|} {$N$ } & \multicolumn{1}{|c|} {$\rho$ }\\ \hline	     %&   (B)   &   (B)    & (B)     \\ \hline 
	     3 refinements   & 4.2967  & 16.7264   & 0.2923  \\ \hline 
	     4 refinements   & 4.6563  & 17.9732   & 0.3136  \\ \hline 
	     5 refinements   & 4.5156  & 26.6800   & 0.4568  \\ \hline 
	\end{tabular}
\end{center}
\caption{Study of the mesh independence of the solver for 3 different levels of refinement.}
\label{mesh}
\end{table}
This table also shows the average value of $s_i^{max}$, named $s^{max}$, for the three refinement levels considered.
From Table \ref{mesh}, it can be seen that the values of $s^{max}$ are quite steady, while the values of $N$ and 
$\rho$ slightly increase as the number of uniform refinements grows.
This increase in the values of $\rho$ indicates that the solver is not completely independent from the mesh considered. A weak dependence on the mesh was expected, since the geometry considered is quite complex, so the mesh constructed on such a geometry is non-trivial and may have low quality elements from the very beginning.   

Table \ref{process} shows the values of $s^{max}$, $N$ and $\rho$ for the case of 4 uniform mesh refinements solved using 18, 36 and 72
processes.
\begin {table}[!h] %[!tpb]
\setlength\tabcolsep{5.5pt} % default value: 6pt
\begin{center}
	\begin{tabular}{|c|c|c|c|c|} \hline	
	\multirow{2}{*} {}   & \multicolumn{3}{|c|}{Process Independence of the Solver} \\ \cline {2-4}
	& \multicolumn{1}{|c|} {$s^{max}$ } & \multicolumn{1}{|c|} {$N$ } & \multicolumn{1}{|c|} {$\rho$ }\\ \hline	     %&   (B)   &   (B)    & (B)     \\ \hline 
	     18 processes   & 4.4688  &  17.5403  & 0.3064  \\ \hline 
	     36 processes   & 4.6563  &  17.9732  & 0.3136  \\ \hline 
	     72 processes   & 4.6563  &  20.6409  & 0.3587  \\ \hline 
	\end{tabular}
\end{center}
\caption{Study of the process independence of the solver for 18, 36 and 72 processes.}
\label{process}
\end{table}
From Table \ref{process}, we see that the performances of the solver are very similar for the three cases considered.
The average number of non-linear iterations $s^{max}$ is almost the same for all the cases, and the value of $N$ increases only three units when using $72$ processes instead of $18$. Finally, the value of $\rho$ grows only by a factor of $10^{-2}$ 
indicating that the solver is process independent.  

\begin{figure}[!htb]
\begin{center}
\begin{minipage}[t]{0.7\linewidth}
\centering
\includegraphics[width=\textwidth]{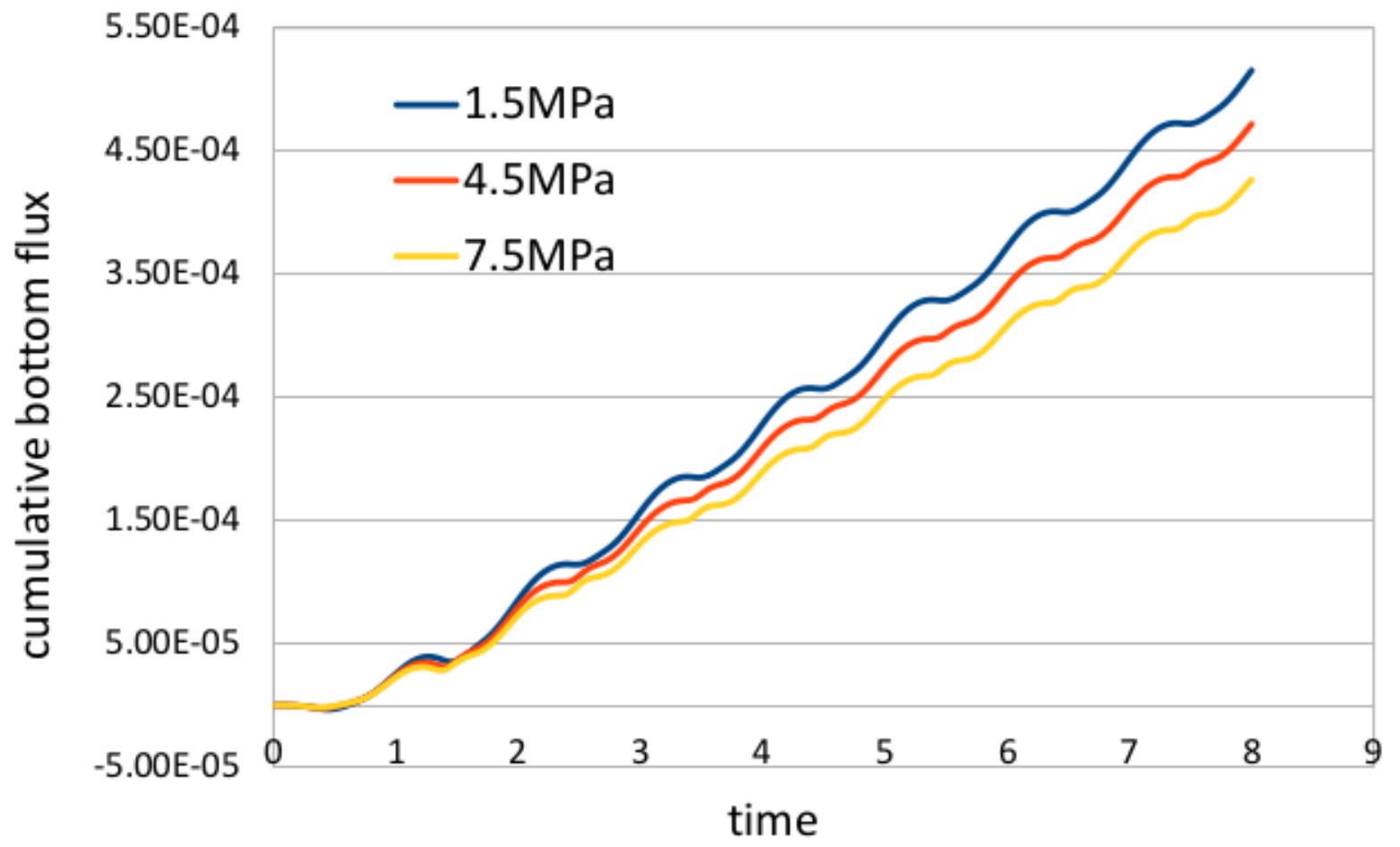} 
 \caption{Cumulative total flux over time at the bottom boundary $B1$.}\label{qBottomIntegral}
 \vspace{0.5cm}
\end{minipage} 
\hspace{.1in}
\begin{minipage}[t]{0.7\linewidth}
\centering
\includegraphics[width=\textwidth]{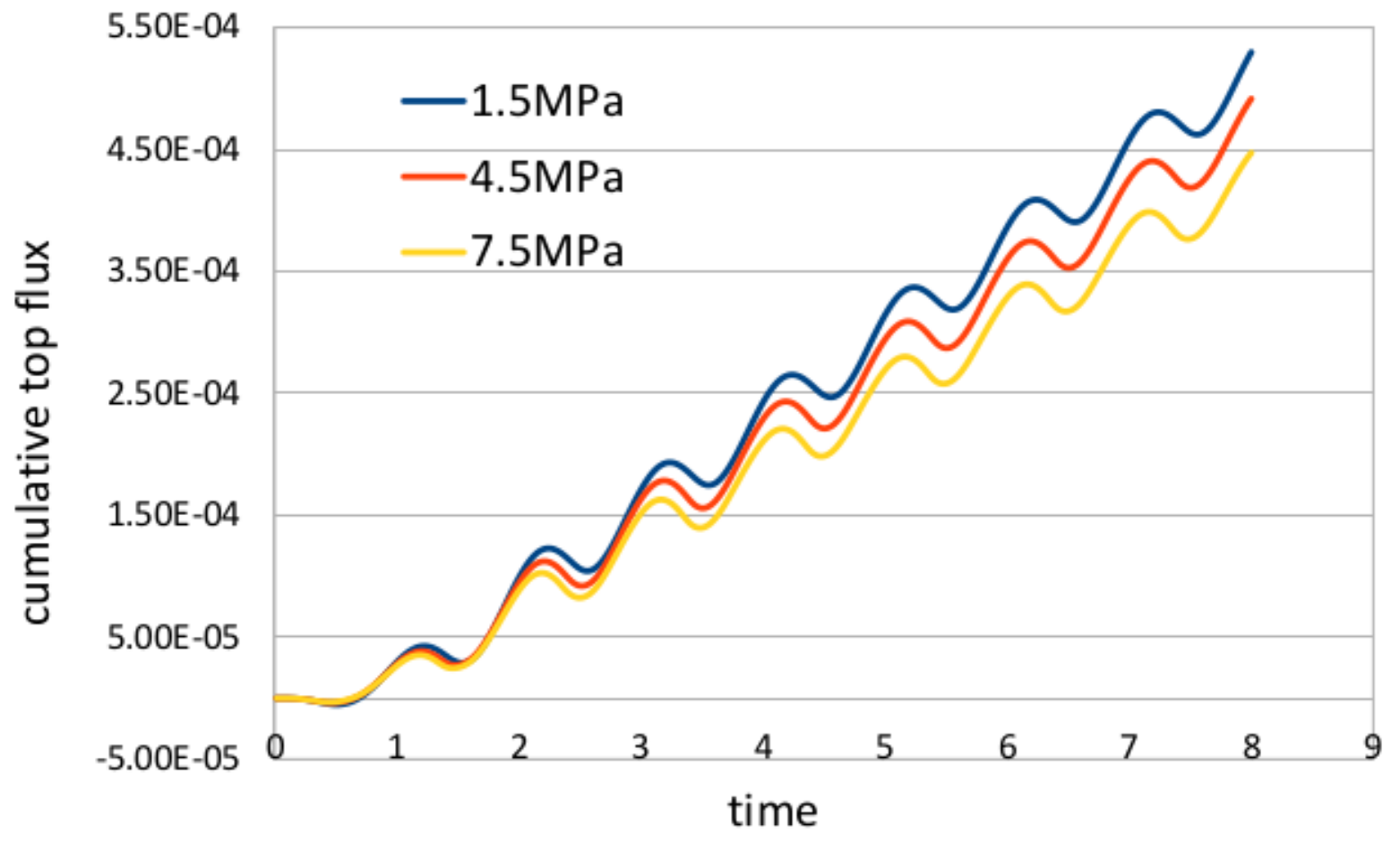}   
\caption{Cumulative total flux over time at the top boundary $B2$.}\label{qTopIntegral}
\end{minipage}
\end{center}
\end{figure}

\subsection{Test 2: sensitivity of valve dynamics to the elastic modulus
of the valve leaflet.}
\begin{figure}[!htb]
\begin{center}
\begin{minipage}[t]{0.7\linewidth}
\centering
\includegraphics[width=\textwidth]{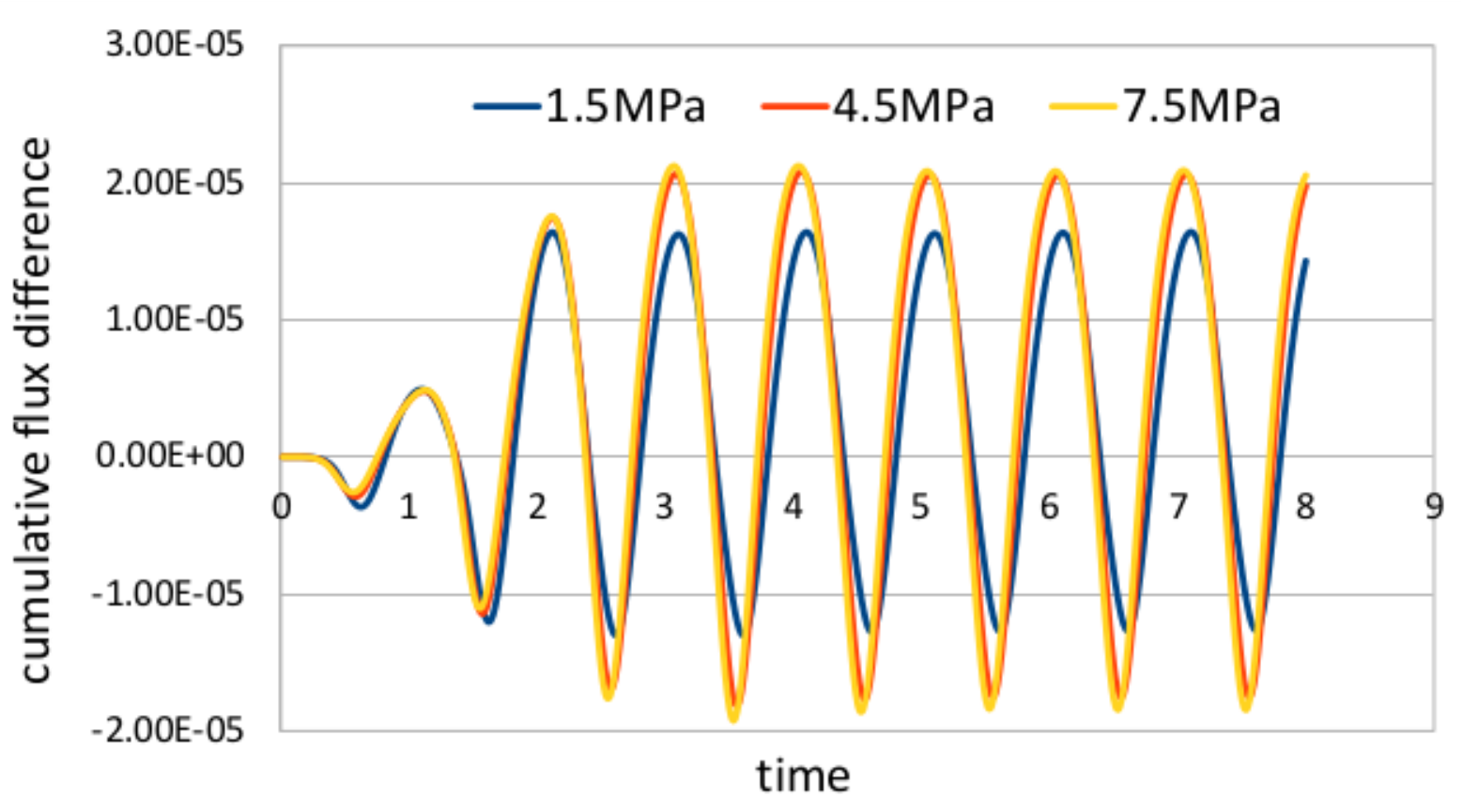}
 \caption{Difference between the cumulative bottom and top fluxes over time.}\label{difference}
 \vspace{0.5cm}
\end{minipage} 
\hspace{.1in}
\begin{minipage}[t]{0.7\linewidth}
\centering
\includegraphics[width=\textwidth]{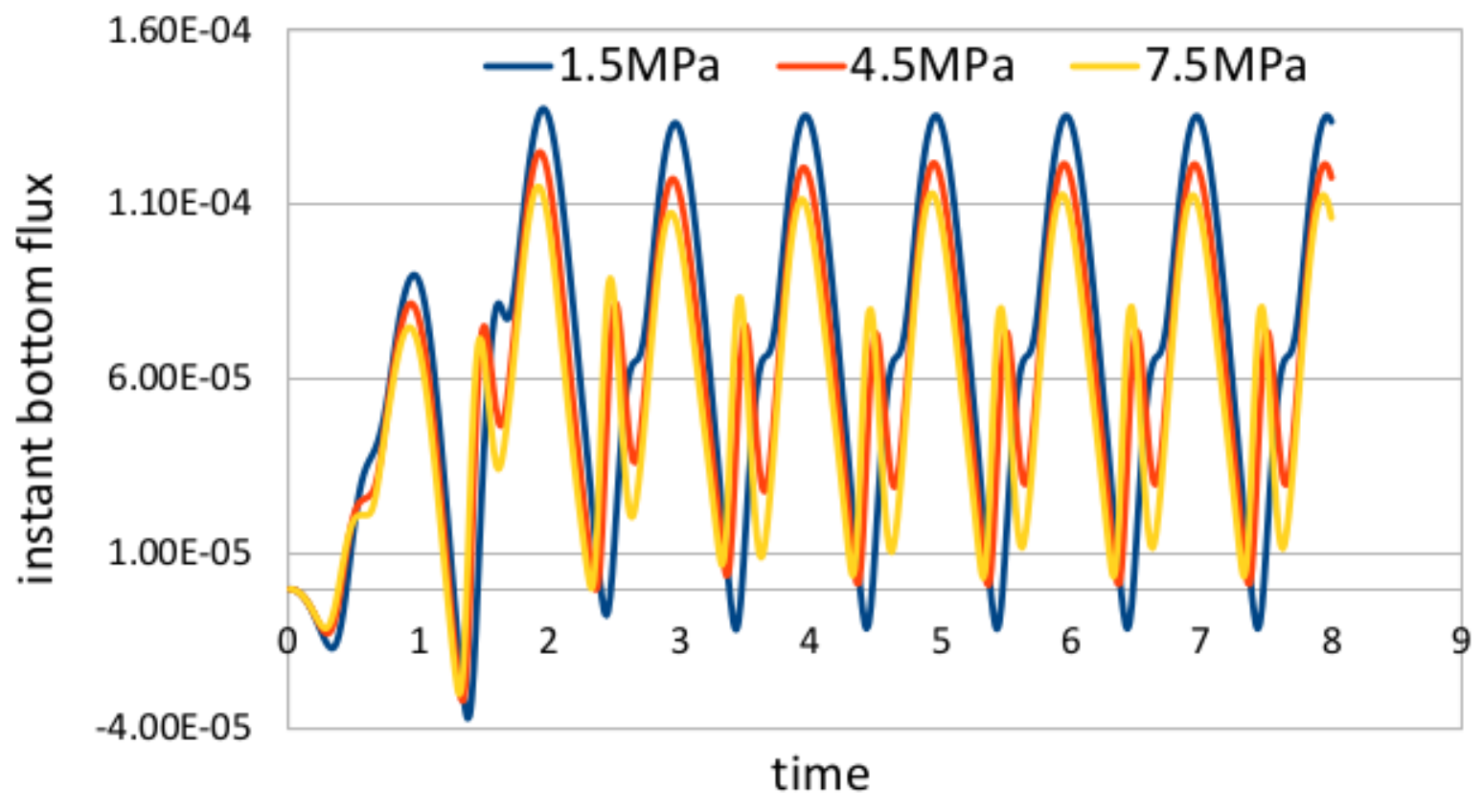}
 \caption{Instant flux at the bottom boundary $B1$.}\label{instantFluxBottom}
\end{minipage}
\end{center}
\end{figure}

Test 2 assesses the sensitivity of valve dynamics to the elastic modulus
of the valve leaflet, independent of vein wall properties. Three analyses are undertaken with different
values for the elastic modulus of the valve leaflet. For the vein wall, a Young's modulus of $260$ MPa is considered, while values of $1.5$, $4.5$ and $7.5$ MPa are specified for the valve leaflet. As pointed out in section \ref{test1}, these parameters do not represent realistic biological values, which is understandable since we are performing 2D simulations that require much stiffer structure properties in order to have reasonable deformations. Figures \ref{qBottomIntegral}, \ref{qTopIntegral}, \ref{difference} and \ref{instantFluxBottom} show, respectively, the cumulative total flux over time at both the bottom ($B1$) and top ($B2$) boundaries, the difference between the cumulative bottom and top fluxes over time and the instant flux at the bottom boundary $B1$ for the three valve leaflet Young's moduli considered. In all the simulations, a time step of $\dfrac{1}{64}$s is employed.
Total flux (instant flux) at a boundary $Bi$, $i=1, 2$, is defined as \begin{equation}
q(t)= \int_{Bi} \bm{u}(x,t) \cdot \mathbf{n} \; ds,                                                             
\end{equation}
while the cumulative total flux over time is the integral over time of the total flux, namely \begin{equation}
Q(t)= \int_{0}^t q(\tau) \; d\tau.                                                             
\end{equation}

\begin{figure}[!htb]
\begin{center}
\begin{minipage}[t]{0.7\linewidth}
\centering
\includegraphics[width=\textwidth]{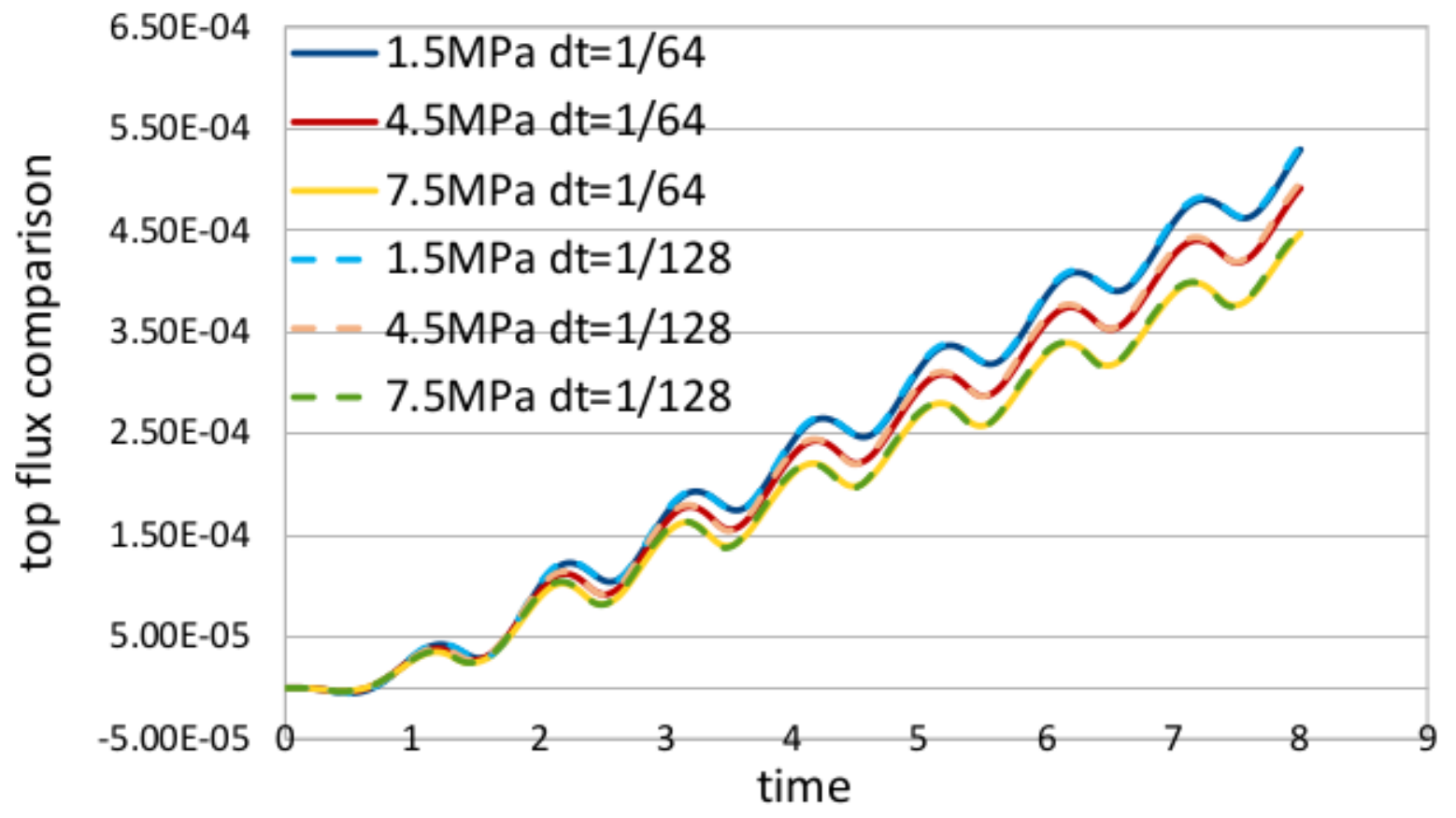}
 \caption{Comparison between the cumulative total flux over time at the top boundary $B2$ using $dt=\frac{1}{64}$s 
 and $dt=\frac{1}{128}$s.}\label{qTopIntegralComparison}
 \vspace{0.5cm}
\end{minipage} 
\hspace{.1in}
\begin{minipage}[t]{0.7\linewidth}
\centering
\includegraphics[width=\textwidth]{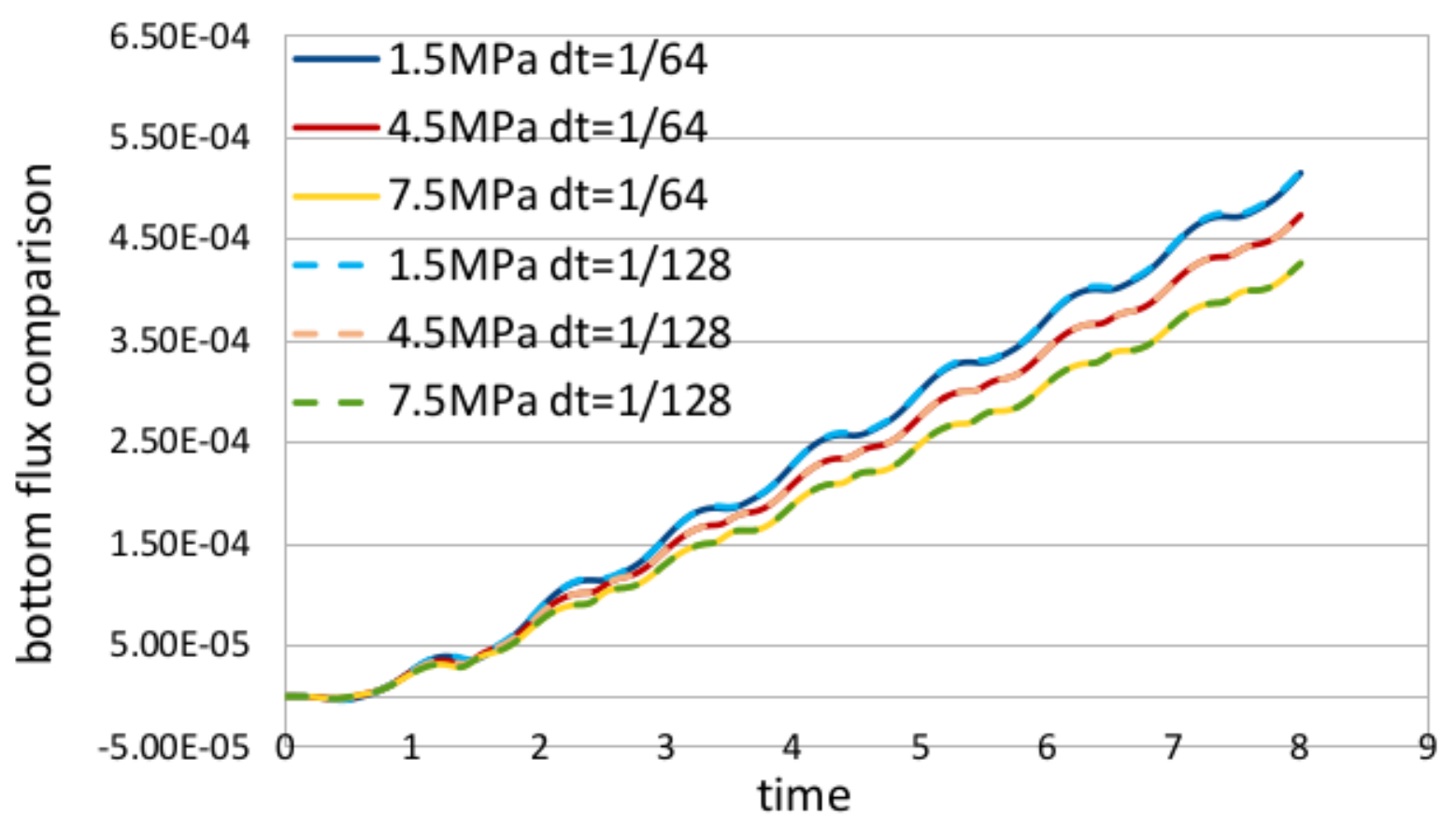}
 \caption{Comparison between the cumulative total flux over time at the bottom boundary $B1$ using $dt=\frac{1}{64}$s 
 and $dt=\frac{1}{128}$s.}\label{qBottomIntegralComparison}
\end{minipage}
\end{center}
\end{figure}

From Figure \ref{difference}, we see that the difference between the cumulative bottom and top fluxes over time always oscillates around $0$ for all three simulations, meaning that the fluxes at the boundaries $B1$ and $B2$ of the vein correctly balance each other over time. Looking closely, it can be seen that the oscillations grow when increasing the leaflet Young's modulus. Figures \ref{qBottomIntegral} and \ref{qTopIntegral} show the cumulative total flux over time  at the boundaries $B1$ and $B2$ of the vein, respectively. In both graphs, the smallest values are attained by the curve corresponding to $7.5$ MPa, meaning by the curve corresponding to the highest Young's modulus value considered. This result indicates that the stiffer the leaflet is, the less blood can flow from the lower to the upper part of the body. In other words, more efficiency is gained with a flexible valve. It is also important to notice that all the curves in both graphs have an oscillating trend. When a curve is locally decreasing, it means that backflow is occurring. One of the main contributions of venous valves is to prevent backflow, but as Figures \ref{qBottomIntegral} and \ref{qTopIntegral} show, this phenomena cannot be completely removed. To have a better understanding of backflow, the instant flux at the bottom ($B1$) of the vein geometry is shown in Figure \ref{instantFluxBottom}. From the graph, it can be clearly seen that backflow is present in all the three cases considered and, the more flexible the leaflet is, the more backflow increases. This growth in the backflow when the valve is more flexible does not contradict the results from Figures \ref{qBottomIntegral} and \ref{qTopIntegral}, because with a flexible valve more blood is pushed from the lower to the upper part of the body, causing higher chances of backflow. 
%Therefore, the right balance between high efficiency and small backflow has to be found     
%\sara{FINIRE DISCORSO SU BACKFLOW.}  

\begin{figure}[!htb]
\begin{center}
 (a) \includegraphics[scale=0.41]{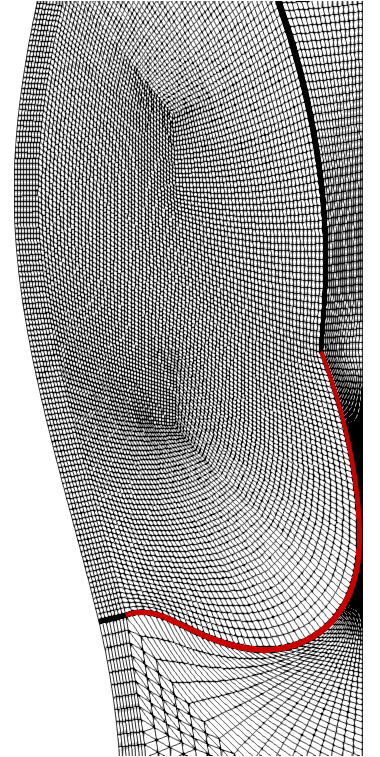} \hspace{0.08cm}
 (b) \includegraphics[scale=0.41]{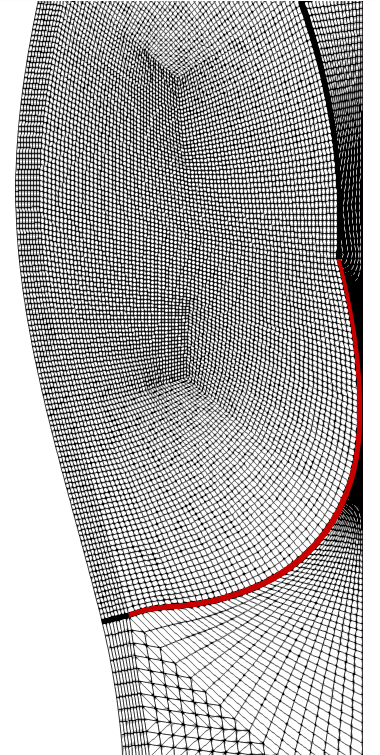} \hspace{0.08cm}
 (c) \includegraphics[scale=0.41]{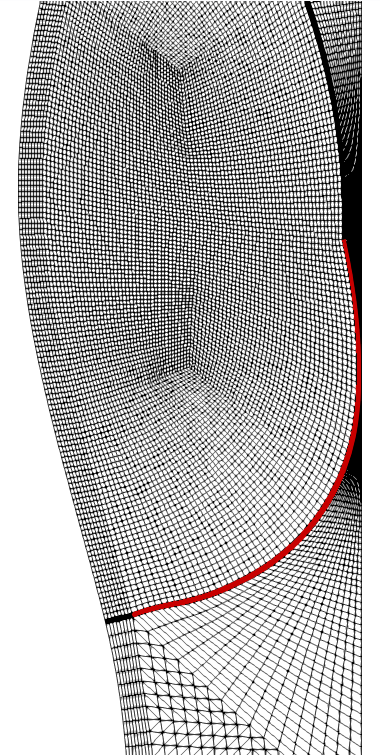}
 \caption{Maximum deformation in the closing phase for (a) $E=1.5$, (b) $E=4.5$ and (c) $E=7.5$ MPa.}\label{max_deforation_closing}
\end{center}
\end{figure}
\begin{figure}[!htb]
\begin{center}
 (a) \includegraphics[scale=0.41]{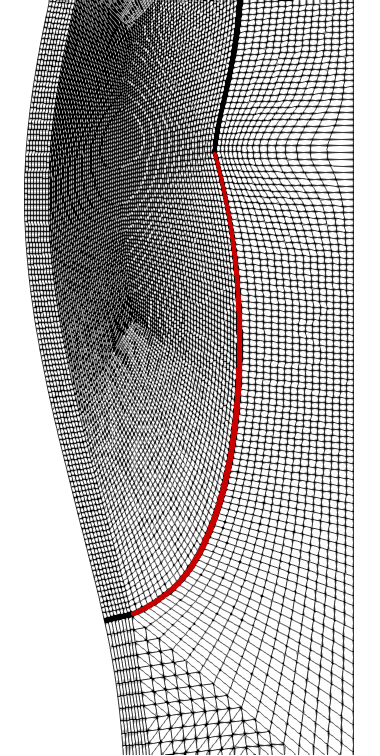} \hspace{0.08cm}
 (b) \includegraphics[scale=0.41]{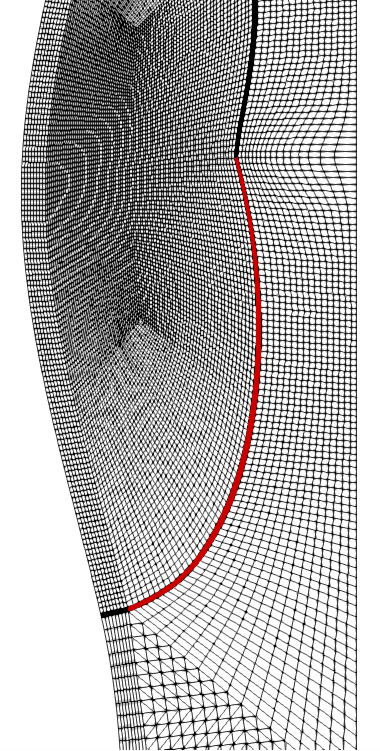} \hspace{0.08cm}
 (c) \includegraphics[scale=0.41]{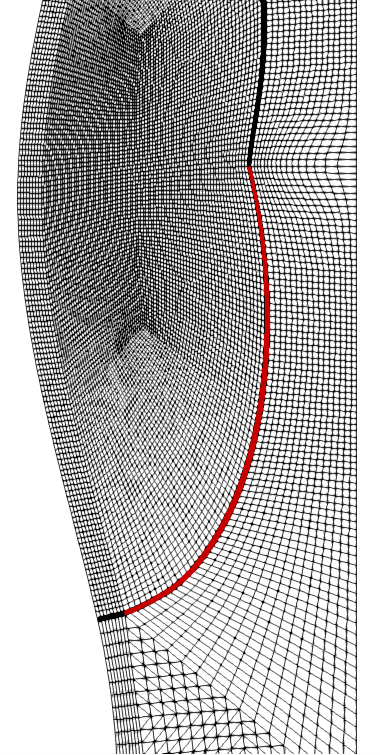}
 \caption{Maximum deformation in the opening phase for (a) $E=1.5$, (b) $E=4.5$ and (c) $E=7.5$ MPa.}\label{max_deforation_opening}
\end{center}
\end{figure}
\begin{figure}[!htb]
\begin{center}
 \includegraphics[scale=0.42]{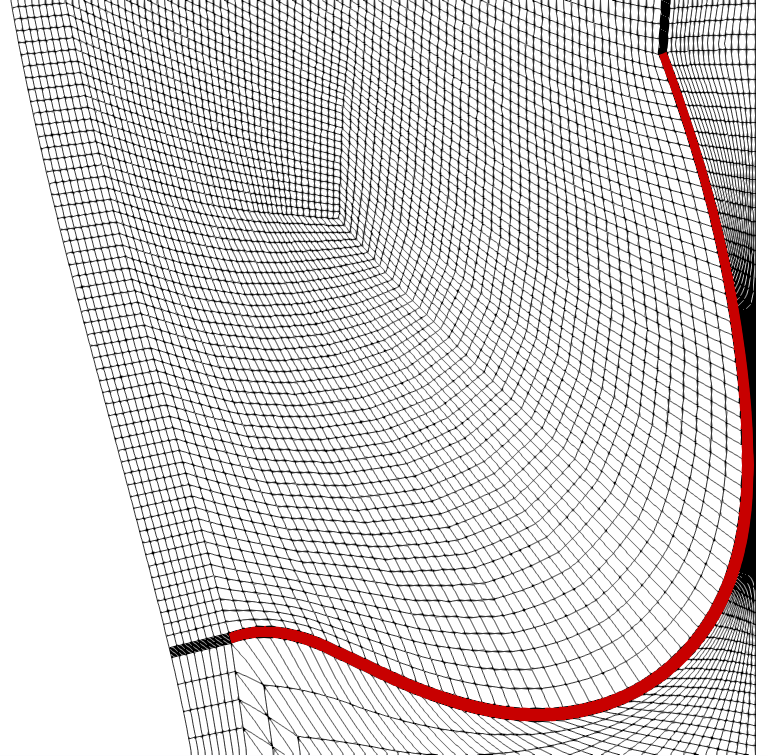}
 \caption{Zoom of the maximum deformation in the opening phase for the $E=1.5$ MPa case.}\label{zoom_max_deforation_opening}
\end{center}
\end{figure}

To study the independence of the solution on the time step, Figures \ref{qTopIntegralComparison} and \ref{qBottomIntegralComparison} compare the cumulative total fluxes over time at the boundaries $B1$ and $B2$ shown before with those obtained using a time step of $\frac{1}{128}$s. The results obtained using a time step of $\frac{1}{128}$s completely overlap the results obtained with a smaller time step of $\frac{1}{64}$s, indicating that there is no dependance of the solution on the time step.

To conclude this section, we show the maximum deformation of the valve leaflet in both the closing and opening phase.
To clearly distinguish the valve leaflet from the surrounding mesh, the leaflet is depicted in red.
From Figure \ref{max_deforation_closing}, we can see that, with a flexible valve ($1.5$ MPa), the deformation is quite heavy, but we still have stability of the simulation. A zoom of this case is shown in Figure \ref{zoom_max_deforation_opening} to provide a clearer vision of the surrounding mesh. From Figure \ref{zoom_max_deforation_opening}, it can be seen that the leaflet gets very close to the axis of symmetry without touching it, because of the action of the fictitious springs.
In both Figures \ref{max_deforation_closing} and \ref{max_deforation_opening}, the \textit{fluid leaflet} is also visible. 

\begin{figure}[!htb]
\begin{center}
 \includegraphics[scale=0.55]{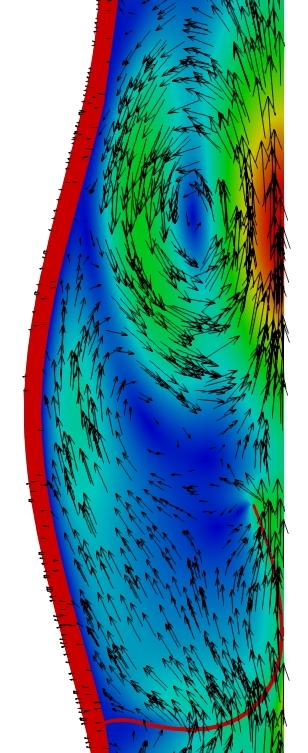} \hspace{1.5cm}
 \includegraphics[scale=0.55]{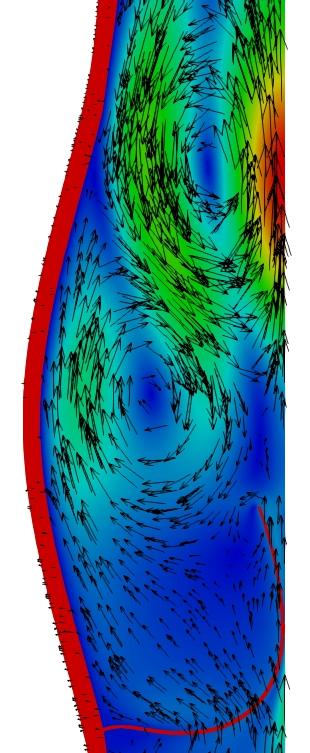}
 \caption{Vortex configurations at time step 100 (left) and 164 (right) for the case $E=1.5$ MPa.}\label{vortexes}
\end{center}
\end{figure}

Finally, Figure \ref{vortexes} portrays the vortex configurations for the case of $E=1.5$ MPa at times $t=1.5625$ s (time step $=100$) and $t=2.5625$ s (time step $=164$). The difference in time between the two snapshots is exactly one valve period.
As the two snapshots show, vortexes occur during the closing phase of a valve cycle. In both configurations, two vortexes are present, but the vortex closer to the venous sinus exhibits different intensities, i.e. is much weaker at time $t=1.5625$ s (time step $=100$) than at time $t=2.5625$ s (time step $=164$).

\section{Conclusions}\label{fine}
Fluid-structure interaction simulations of venous valves are a challenging problem since the large structural displacements
of the valve leaflets may lead to mesh deteriorations and entanglements, causing instabilities of the solver and, consequently, the
numerical solution to diverge. In this paper, we present a monolithic ALE scheme
for FSI simulations of venous valves designed to solve these instabilities. The proposed method is based on three main features: a staggered in time mesh velocity in the discretization scheme to improve computational stability; a scaling factor that measures the distance of a fluid element from the valve leaflets, to guarantee that there are no mesh entanglements in the fluid domain; and fictitious springs to model the contact force between closing valve leaflets. To further improve stability, a Streamline Upwind
Petrov Galerkin stabilization is added to the momentum equation.
The Newton-Krylov solver employed is described, where
we consider the use of geometric multigrid preconditioners. We describe the structure of the
geometric multigrid operators, for which modified Richardson smoothers are chosen, preconditioned
by an additive Schwarz algorithm of overlapping restricted type. 

We perform several 2D tests to assess the proposed method. The tests about the solver performances show process independence of the solver and a weak dependance of the solver on the mesh. 
To assess the sensitivity of valve dynamics to the elastic modulus of the valve
leaflets, independent of vein wall properties, three analyses are performed with different values
for the elastic modulus of the valve leaflet, namely $1.5$, $4.5$ and $7.5$ MPa. These tests show that 
cumulative total fluxes over time and backflow are lower when increasing the valve Young's modulus. 
Thus, with a flexible valve, more blood can flow from the lower to the upper part of the body, but there could be a higher backflow. This result shows that valve dynamics is very sensible to the elastic modulus of the valve leaflets, and that the flexibility of the leaflets plays a central role in the valve mechanism.

Future work will consist in performing 3D simulations of venous valves and investigating a solver that could provide better computational performances.  

\section*{Acknowledgments}
This work was supported by the National Science Foundation grant DMS-1412796.

\bibliographystyle{wileyj}
\bibliography{particle}

\end{document}